\documentclass[reqno,12pt]{amsart}

\usepackage[spanish,english]{babel}
\usepackage[dvips]{graphicx}
\usepackage{amsmath,multirow}
\usepackage{pstricks}
\usepackage[pdflatex]{epsfig}
\usepackage[utf8]{inputenc}
\usepackage[T1]{fontenc}
\newcommand{\minmod}{\operatorname*{\mathrm{mm}}}

\allowdisplaybreaks
\numberwithin{equation}{section}

\title[High order numerical schemes for non-local conservation laws]{High order numerical schemes for one-dimension non-local conservation laws}
\author[Chalons]{Christophe Chalons$^{\mathrm{a}}$} 
\author[Goatin]{Paola Goatin$^{\mathrm{b}}$} 
\author[Villada]{Luis M. Villada$^{\mathrm{c}}$}
 \thanks{$^{\mathrm{a}}$Laboratoire de Math\'ematiques de Versailles, UMR~8100, Universit\'e de
Versailles Saint-Quentin-en-Yvelines, UFR des Sciences, b\^{a}timent
Fermat, 45 avenue des Etats-Unis, 78035 Versailles cedex,
France. E-mail: {\tt christophe.chalons@uvsq.fr}}
\thanks{$^{\mathrm{b}}$ Inria Sophia Antipolis - M\'editerran\'ee, Universit\'e C\^ote d'Azur, Inria, CNRS, LJAD, 2004, route des Lucioles - BP 93 06902 Sophia Antipolis Cedex, France  E-mail: {\tt paola.goatin@inria.fr}}
\thanks{$^{\mathrm{c}}$GIMNAP-Departamento de Matem\'{a}ticas, Universidad del B\'{\i}o-B\'{\i}o,   Casilla 5-C,  Concepci\'{o}n Chile
and  CI$^{\mathrm{2}}$MA, Universidad de Concepci\'{o}n, Casilla 160-C, Concepci\'{o}n, Chile.  E-Mail:    {\tt lvillada@ubiobio.cl}}

\subjclass[2000]{35L65,90B20,76T20,65M08,65M60.}
 \keywords{Scalar conservation laws, non-local flux, sedimentation models, traffic flow model, Galerkin discontinuous schemes, 
 finite volume schemes.}
\date{\today}

\begin{document}

\begin{abstract}
This paper focuses on the numerical approximation of the solutions of non-local conservation laws in one space dimension. 
These equations are motivated by two distinct applications, namely a traffic flow model in which the mean velocity depends on a weighted mean of the downstream traffic density, and a  sedimentation model where either the solid phase velocity or the solid-fluid relative velocity depends on the concentration in a neighborhood. In both models, the velocity is a function of a convolution product between the unknown and a kernel function with compact support. It turns out that the solutions of such equations may exhibit oscillations that are very difficult to approximate using classical first-order numerical schemes. We propose to design Discontinuous Galerkin (DG) schemes and Finite Volume WENO (FV-WENO) schemes to obtain high-order approximations. As we will see, the DG schemes give the best numerical results but their CFL condition is very restrictive. On the contrary, FV-WENO schemes can be used with larger time steps. 
 We will see that the evaluation of the convolution terms necessitates the use of quadratic polynomials reconstructions in each cell in order to obtain the high-order accuracy with the FV-WENO approach. Simulations using DG and FV-WENO schemes are presented for both applications.
\end{abstract}
\maketitle

\section{Introduction}
This paper is concerned with the design of numerical schemes for the one-dimensional Cauchy problem for non-local scalar conservation law of the form
\begin{equation}\label{nl_claw}
\begin{cases}
\rho_t+(f(\rho)V(\rho*\omega_\eta ))_x=0, & x\in \mathbb{R},\quad t>0, \\  %
\rho(x,0)=\rho_0(x) & x\in  \mathbb{R},
\end{cases}
\end{equation}
where the unknown density $\rho$ depends on the space variable $x$ and the time variable $t$, $\rho \to V(\rho)$ is a given velocity function, and  $\rho \to g(\rho) = f(\rho)V(\rho)$ is the usual flux function for the corresponding local scalar conservation law. 
Here, (\ref{nl_claw}) is non-local in the sense that the velocity function $V$ is evaluated on a ``neighborhood'' of $x\in\mathbb{R}$ defined by the convolution of the density $\rho$ and a given kernel function $\omega_\eta$ with compact support. In this paper, we are especially interested in two specific forms of (\ref{nl_claw}) which naturally arise in traffic flow modelling \cite{BG_2016, GS_2016} and sedimentation problems \cite{BBKT_2011}. They are given as follows. \\  
\ \\
{\it A non-local traffic flow model.} In this context, we follow \cite{BG_2016, GS_2016} and consider (\ref{nl_claw}) as an extension of the classical Lighthill-Whitham-Richards traffic flow model, in which the mean velocity is assumed to be 
a non-increasing function of the downstream traffic density and where the flux function is given by
\begin{equation}\label{LWRnl}
f(\rho)=\rho,\quad
V(\rho)= \textcolor{black}{1-\rho}, \quad  (\rho*\omega_\eta)(x,t)=\int_x^{x+\eta}\omega_\eta(y-x)\rho(y,t)dy.
\end{equation}
Four different nonnegative kernel functions will be considered in the numerical section, namely $\omega_\eta(x) = 1/\eta$, 
$\omega_\eta(x) = 2(\eta-x)/\eta^2$,
$\omega_\eta(x) = 3(\eta^2-x^2)/(2 \eta^3)$ and \textcolor{black}{$\omega_\eta(x) =  2x/\eta^2$} with support on $[0,\eta]$ for a given value of the real number $\eta>0$. Notice that the well-posedness of this model together with the design of a first and a second order FV approximation have been considered in \cite{BG_2016, GS_2016}. \\
\ \\
{\it A non-local sedimentation model.} Following \cite{BBKT_2011}, we consider under idealized assumptions that equation \eqref{nl_claw}  represents a one-dimensional model for the sedimentation 
of small equal-sized spherical solid particles dispersed in a viscous fluid, where the local solids column fraction $\rho=\rho(x,t)$ is a function of depth $x$
and time $t$. In this context, the flux function is given by 
\begin{equation}\label{sednl}
f(\rho)=\rho(1-\rho)^\alpha,\quad
V(\rho)=(1-\rho)^n, \quad (\rho*\omega_\eta)(x,t)=\int_{-2\eta}^{2\eta}\omega_\eta(y)\rho(x+y,t)dy,
\end{equation}
where $n \geq 1$ and the parameter $\alpha$ satisfies $\alpha=0$ or $\alpha\geq1$. The function $V$ is the so-called hindered settling factor and the convolution term $\rho*\omega_\eta$ is defined by a symmetric, nonnegative, and piecewise smooth kernel function $\omega_\eta$ with 
support on $[-2\eta,2\eta]$ for a parameter $\eta>0$ and $\int_{\mathbb{R}}\omega_\eta(x)dx=1$. More precisely, the authors define in \cite{BBKT_2011} a truncated parabola $K$ by 
$$
K(x)=\frac38\left(1-\frac{x^2}{4}\right) \, \text{ for } |x|<2,\qquad K(x)=0\,\,\text{ otherwise, }
$$
and set
 \begin{equation}\label{ker_sed}
 \omega_\eta(x):=\eta^{-1}K(\eta^{-1}x).
 \end{equation}

Conservation laws with non-local terms arise in a variety of physical and engineering applications: besides the above cited ones, we mention models for granular flows \cite{AmadoriShen2012}, production networks \cite{Keimer2015}, conveyor belts \cite{Gottlich2014}, weakly coupled oscillators \cite{AHP2016}, laser cutting processes \cite{ColomboMarcellini2015}, biological applications like structured populations dynamics \cite{Perthame_book2007}, crowd dynamics \cite{Carrillo2016, ColomboGaravelloMercier2012} or more general problems like gradient constrained equations \cite{Amorim2012}. \\
While several analytical results on non-local conservation laws can be found in the recent literature (we refer to \cite{Amorim2015}
for scalar equations in one space dimension, \cite{AmbrosioGangbo2008, ColomboHertyMercier2011, PiccoliRossi2013} for scalar equations in several space dimensions and \cite{ACG2015,ColomboMercier2012,CrippaMercier2012} for multi-dimensional systems of conservation laws),
few specific numerical methods have been developed up to now.
Finite volume numerical methods have been studied in \cite{ACG2015, GS_2016, KurganovPolizzi2009, PiccoliTosin2011}. At this regard, it is important to notice that the lack of Riemann solvers for non-local equations limits strongly the choice of the scheme. At the best of our knowledge, two main approaches have been proposed in the literature to treat non-local problems: first and second order central schemes like 
Lax-Friedrichs or Nassyau-Tadmor \cite{ACG2015, Amorim2015, BBKT_2011, GS_2016, KurganovPolizzi2009} 
and Discontinuous Galerkin (DG) methods \cite{GottlichSchindler2015}. In particular, the comparative study presented in \cite{GottlichSchindler2015}
on a specific model for material flow in two space dimensions, involving density gradient convolutions,
encourages the use of DG schemes for their versatility and lower computational cost, but further investigations are needed in this direction. 
Besides that, the computational cost induced by the presence of non-local terms, requiring the computation of quadrature formulas at each time step,
motivate the development of high order algorithms.

The aim of the present article is to conduct a comparison study on high order schemes for a class of non-local scalar equations in one space dimension,
focusing on equations of type \eqref{nl_claw}.
In Section \ref{review} we review DG and FV-WENO schemes for classical conservation laws.
These schemes will be extended to the non-local case in Sections \ref{sec:dg_nl} and \ref{sec:fv_nl}. 
Section \ref{sec:tests} is devoted to numerical tests.

%
%

\section{A review of Discontinuous Galerkin and Finite Volume WENO  schemes for local conservation laws} \label{review}
The aim of this section is to introduce some notations and to briefly review the DG and FV-WENO numerical schemes to solve the classical {\it local} nonlinear conservation law 
\begin{equation}\label{claw}
\begin{cases}
\rho_t+g(\rho)_x=0, & x\in \mathbb{R},\quad t>0, \\
\rho(x,0)=\rho_0(x), & x\in \mathbb{R}.
\end{cases}
\end{equation}

We first consider $\{I_j\}_{j\in\mathbb{Z}}$ a partition of $\mathbb{R}$. The points $x_j$ will represent the centers of the cells $I_j=[x_{j-\frac12},x_{j+\frac12}]$, and the cell sizes will be denoted by 
$\Delta x_j=x_{j+\frac12}-x_{j-\frac12}$. The largest cell size is $h=\sup_j \Delta x_j$. Note that, in practice, we will consider a constant space step so that we will have $h=\Delta x$.

\subsection{The Discontinuous Galerkin approach}

In this approach, we look for approximate
solutions in the function space of discontinuous polynomials
$$
V_h:=V_h^k=\{v:v|_{I_j}\in P^k(I_j)\},
$$
where $P^k(I_j)$ denotes the space of polynomials of degree at most $k$ in the element $I_j$.
The approximate solutions are sought under the form
$$
\rho^h(x,t) = \sum_{l=0}^k c_j^{(l)}(t) v_j^{(l)}(x), \quad v_j^{(l)}(x) = v^{(l)}(\zeta_j(x)),
$$ 
where $c_j^{(l)}$ are the degrees of freedom 
in the element $I_j$.
The subset $\{v_j^{(l)}\}_{l=0,...,k}$ 
constitutes a basis of $P^k(I_j)$ and 
in this work we will take Legendre polynomials as a local orthogonal basis of $P^k(I_j)$, namely
$$v^{(0)}(\zeta_j)=1,\qquad v^{(1)}(\zeta_j)=\zeta_j ,\qquad   v^{(2)}(\zeta_j)=\frac12\left( 3 \zeta_j^2-1\right),\dots, \quad \zeta_j:=\zeta_j(x)=\frac{x-x_j}{{\Delta x}/2},$$
see for instance \textcolor{black}{\cite{CS_2001,QS_2005}}. \\
Multiplying \eqref{claw} by $v_h\in V_h$ and integrating over $I_j$ gives 
\begin{equation}
\int_{I_{j}}\rho_{t}v_{h}dx-\int_{I_{j}}g(\rho)v_{h,x}dx+g(\rho(\cdot,t))v_{h}(\cdot)\bigl\lvert_{x_{j-\frac12}}^{x_{j+\frac12}}=0,
\end{equation}
and {the semi-discrete DG formulation thus consists} in looking for $\rho^{h}\in V_{h}$, such that for all $v_{h}\in V_{h}$ and all $j$, 
\begin{equation}\label{semi-discrete}
\int_{I_{j}}\rho^h_{t}v_{h}dx-\int_{I_{j}}g(\rho^h)v_{h,x}dx+\hat{g}_{j+\frac12}v^{-}_{h}(x_{j+\frac12})-\hat{g}_{j-\frac12}v^{+}_{h}(x_{j-\frac12})=0 ,
\end{equation}
where $\hat{g}_{j+\frac12}=\hat{g}(\rho^{h,-}_{j+\frac12},\rho^{h,+}_{j+\frac12})$ is a consistent, monotone and Lipschitz continuous numerical flux function. In particular, 
we will choose to use the Lax-Friedrichs flux
$$\hat{g}(a,b):=\frac{g(a)+g(b)}{2}+\alpha\frac{a-b}{2},\quad\quad \alpha=\max_u |g'(u)|.$$
Let us now observe that if $v_h$ is the $l$-th Legendre polynomial, we have $v^{+}_{h}(x_{j-\frac12})=v^{(l)}(\zeta_{j}(x_{j-\frac12}))=(-1)^{l}$, and 
$v^{+}_{h}(x_{j+\frac12})=v^{(l)}(\zeta_{j}(x_{j+\frac12}))=1$, $\forall j\,, l=0,1,\dots,k$. Therefore, 
replacing $\rho^h(x,t)$ by $\rho_j^{h}(x,t)$ and $v_h(x)$ by $v^{(l)}(\zeta_j(x))$ in \eqref{semi-discrete}, the degrees of freedom $c^{(l)}_{j}(t)$ satisfy the differential equations 
\begin{equation}\label{ode}
\frac{d}{dt}c^{(l)}_{j}(t)+\frac{1}{a_{l}} \left(-\int_{I_{j}}g(\rho_j^h)\frac{d}{dx}v^{(l)}(\zeta_j(x))dx+\hat{g}(\rho^{h,-}_{j+\frac12},\rho^{h,+}_{j+\frac12})-
(-1)^{l}\hat{g}(\rho^{h,-}_{j-\frac12},\rho^{h,+}_{j-\frac12})\right)=0,
\end{equation}
with $${a_{l}}=\int_{I_{j}}(v^{(l)}(\zeta_j(x)))^2dx=\frac{\Delta x}{2l+1},\qquad l=0,1,\dots,k.$$ 
On the other hand, the initial condition can be obtained using the $L^2$-projection of $\rho_0(x)$, namely
$$c^{(l)}_{j}(0)=\frac{1}{a_{l}}\int_{I_{j}}\rho_0(x)v^{(l)}(\zeta_j(x))dx=\frac{2l+1}{2}\int_{-1}^1\rho_0\left(\frac{\Delta x}{2}y+x_j\right)v^{(l)}(y)dy,\qquad l=0,\dots,k.$$
The integral terms in \eqref{ode} can be computed exactly or using a high-order quadrature technique after a suitable change of variable, namely
 $$\int_{I_{j}}g(\rho_j^h)\frac{d}{dx}v^{(l)}(\zeta_j(x))dx=\int_{-1}^1g\left(\rho_j^h\left(\frac{\Delta x}{2}y+x_j,t\right)\right)(v^{(l)})'(y)dy.$$ 
In this work,  we will consider a Gauss-Legendre quadrature with $N_G=5$  nodes for integrals on $[-1,1]$
$$
\int_{-1}^{1}g(y)dy=\sum_{e=1}^{N_G}w_eg(y_e),
$$ 
where $y_e$ are the Gauss-Legendre quadrature points such that the quadrature formula is exact for polynomials of degree until $2N_G-1=9$ \textcolor{black}{\cite{abramowitz1966handbook}}. \\

\textcolor{black}{
The semi-discrete scheme \eqref{ode} can be written under the usual form  
$$\frac{d}{dt}C(t)=	\mathcal{L}(C(t)),$$ 
where $\mathcal{L}$ is the spatial discretization operator defined by (\ref{ode}). In this work, we will consider a time-discretisation based on the following
 total variation diminishing (TVD) third-order Runge-Kutta  method \cite{shu1988efficient},
\begin{eqnarray}\label{rk3}
C^{(1)}&=&C^{n}+\Delta t \mathcal{L}(C^{n}),\notag\\
C^{(2)}&=&\frac34C^{n}+\frac14C^{(1)}+ \frac14\Delta t \mathcal{L}(C^{(1)}),\\
C^{n+1}&=&\frac13C^{n}+\frac23C^{(2)}+ \frac23\Delta t \mathcal{L}(C^{(2)}).\notag
\end{eqnarray}
Other TVD or strong stability preserving SSP time discretization can be also used \cite{gottlieb2009high}.  The CFL condition is 
$$\frac{\Delta t}{\Delta x}\max_\rho|g'(\rho)|\leq \mathrm{C}_{\mathrm{CFL}}= \frac{1}{2k+1},$$
where $k$ is the degree of the polynomial, see \textcolor{black}{\cite{CS_2001}}. The scheme 
\eqref{ode} and \eqref{rk3} will be denoted RKDG. }

\subsection{Generalized slope limiter}
It is well-known that RKDG schemes like the one proposed above may oscillate when sharp discontinuities are present in the solution. In order to control these instabilities, a common
strategy is to use a limiting procedure. We will consider the so-called generalized slope limiters proposed in \cite{CS_2001}. 
With this in mind and $\rho_j^h(x)=\sum_{l=0}^{k}c^{(l)}_{j}v^{(l)}(\zeta_j(x))\in P^k(I_j)$, we first set
$$u_{j+\frac12}^{-}:=\bar{\rho}_j+\minmod(\rho_j^h(x_{j+\frac12})-\bar{\rho}_j,\Delta_{+}\bar{\rho}_j,\Delta_{-}\bar{\rho}_{j})
$$
and
$$ u_{j-\frac12}^{+}:=\bar{\rho}_j-\minmod(\bar{\rho}_j-\rho_j^h(x_{j-\frac12}),\Delta_{+}\bar{\rho}_j,\Delta_{-}\bar{\rho}_{j}),$$
where $\bar{\rho}_j$ is the average of $\rho^h$ on $I_j$,   $\Delta_{+}\bar{\rho}_j=\bar{\rho}_{j+1}-\bar{\rho}_j$,  $\Delta_{-}\bar{\rho}_j=\bar{\rho}_j-\bar{\rho}_{j-1}$, and 
where $\minmod$ is given by the minmod function limiter
$$\minmod(a_1,a_2,a_3)=\begin{cases} s\cdot \min_{j} |a_j| & \text{ if } s=sign(a_1)=sign(a_2)=sign(a_3)\\ 0 & \text{ otherwise,}    \end{cases}$$
or by the TVB modified minmod function 
\begin{equation}\label{slopelimiter}
 {\overline{\minmod}}(a_1,a_2,a_2)=\begin{cases} a_1 & \text{ if } |a_1|\leq M_{b} h^2,\\ \minmod(a_1,a_2,a_3) & \text{ otherwise,}    \end{cases}
\end{equation}
where $M_{b} >0$ is a constant. According to \cite{CS_2001, QS_2005}, this constant is proportional to the second-order derivative of the initial condition at local extrema.
Note that if $M_b$ is chosen too small, the scheme is very diffusive, while if $M_b$ is too large, oscillations may appear. \\
The values $u_{j+\frac12}^{-}$ and $u_{j-\frac12}^{+}$ allow to compare the interfacial values of $\rho^h_j(x)$ with respect to its local cell averages. Then, the generalized slope limiter technique consists in replacing $\rho_j^h$ on each cell $I_j$ with $\Lambda\Pi_h$ defined by
\begin{equation*}
\Lambda\Pi_h(\rho_j^h)=\begin{cases}
 \rho_j^h &\text{if } u_{j-\frac12}^{+}=\rho_j^h(x_{j-\frac12})\text{ and } u_{j+\frac12}^{-}=\rho_j^h(x_{j+\frac12}),\\ 
\displaystyle \bar{\rho}_j+\frac{(x-x_j)}{\Delta x/2}\minmod(c^{(1)}_{j},\Delta_{+}\bar{\rho}_j,\Delta_{-}\bar{\rho}_{j}) &\text{otherwise.}  
     \end{cases} 
\end{equation*}
Of course, this generalized slope limiter procedure has to be performed after each inner step of the Runge-Kutta scheme \eqref{rk3}.

\subsection{The Finite Volume WENO approach} \label{sec:fvweno1}
In this section, we solve the nonlinear conservation law \eqref{claw} by using a 
high-order finite volume WENO scheme \textcolor{black}{\cite{Shu_1998,shu1988efficient}}. 
Let us denote by $\bar{\rho}(x_j,t)$ the cell average of the exact solution $\rho(\cdot,t)$ in the cell $I_j$~:
$$
\bar{\rho}(x_j,t):=\frac{1}{\Delta x}\int_{I_j}\rho(x,t)dx.
$$
The unknowns are here the set of all $\{\bar{\rho}_j(t)\}_{j\in\mathbb{Z}}$ which represent approximations of the exact cell averages $\bar{\rho}(x_j,t)$. Integrating \eqref{claw} over $I_j$
we obtain 
\begin{equation*}
\frac{d}{dt}\bar{\rho}(x_j,t)=-\frac{1}{\Delta x}\left( g(\rho(x_{j+\frac12},t))-g(\rho(x_{j-\frac12},t))\right),\qquad \forall j \in \mathbb{Z}.
\end{equation*}
This equation is approximated by the semi-discrete conservative scheme
\begin{equation}\label{fv_semi}
\frac{d}{dt}\bar{\rho}_j(t)=-\frac{1}{\Delta x}\left( \hat{g}_{j+\frac12}-\hat{g}_{j-\frac12}\right),\qquad \forall j \in \mathbb{Z},
\end{equation}
where the numerical flux $\hat{g}_{j+\frac12}:=\hat{g}(\rho_{j+\frac12}^l,\rho_{j+\frac12}^r)$ is the Lax-Friedrichs  flux and $\rho_{j+\frac12}^l$ and $\rho_{j+\frac12}^r$ are some left and right high-order WENO reconstructions of $\rho(x_{j+\frac12},t)$ obtained from the cell averages $\{\bar{\rho}_j(t)\}_{j \in \mathbb{Z}}$. 
Let us focus  on the definition of $\rho_{j+\frac12}^l$. 
In order to obtain a $(2k-1)$th-order WENO approximation, we first compute $k$ reconstructed values
$$\hat{\rho}_{j+\frac12}^{(r)}=\sum_{i=0}^{k-1}c_{r}^i\bar{\rho}_{i-r+j},\quad r=0,\dots,k-1,$$
that correspond to $k$ possible stencils $S_r(j)=\{x_{j-r},\dots,x_{j-r+k-1}\}$ for $r=0,\dots,k-1$. The coefficients $c_{r}^i$ are chosen in such a way that each of the $k$ reconstructed values is $k$th-order accurate \cite{Shu_1998}. Then, the $(2k-1)$th-order WENO reconstruction is a convex combination of all these $k$ reconstructed values and defined by
$$
\rho_{j+\frac12}^l=\sum_{r=0}^{k-1}w_{r}\hat{\rho}_{j+\frac12}^{(r)},
$$
where the positive nonlinear weights $w_r$ satisfy $\sum_{r=0}^{k-1}w_r=1$ and are defined by
$$
w_r=\frac{\alpha_r}{ \sum_{s=0}^{k-1}\alpha_s},\quad \alpha_r=\frac{d_r}{(\epsilon+\beta_r)^2}.
$$
Here $d_r$ are the linear weights which yield the $(2k-1)$th-order accuracy, $\beta_r$ are called the ``smoothness indicators'' of the stencil $S_r(j)$, which measure the smoothness of the function $\rho$ in the stencil, and $\epsilon$ is a small parameter used to avoid dividing by zero (typically $\epsilon=10^{-6}$).
The exact form of the smoothness indicators and other details about WENO reconstructions can be found in \cite{Shu_1998}. \\
The reconstruction of $\rho_{j-\frac12}^r$ is obtained in a mirror symmetric fashion with respect to $x_{j-\frac12}$.
\textcolor{black}{
The semi-discrete scheme \eqref{fv_semi} is then integrated in time using the (TVD) third-order Runge-Kutta scheme \eqref{rk3}, under 
the CFL condition 
$$\frac{\Delta t}{\Delta x}\max_\rho|g'(\rho)|\leq \mathrm{C}_{\mathrm{CFL}}<1.$$. }
%
%

\section{Construction of DG schemes for non-local problems}  \label{sec:dg_nl}
We now focus on the non-local equation \eqref{nl_claw}, for which we set $R(x,t):=(\rho*\omega_\eta)(x,t)$. Since  
$\rho^{h}(x,t)|_{I_j}=\sum_{l=0}^{k}c^{(l)}_{j}(t)v^{(l)}(\zeta_j(x))\in P^k(I_j)$, is a weak solution of the non-local problem  \eqref{nl_claw}, the coefficients $c^{(l)}_{j}(t)$ can be calculated by solving the following differential equation,
\begin{equation}\label{DG_nl}
\frac{d}{dt}c^{(l)}_{j}(t)+\frac{1}{a_{l}} \left(-\int_{I_{j}}f(\rho^h) V(R)\frac{d}{dx}v^{(l)}(\zeta_j(x))dx+\check{f}_{j+\frac12}-
(-1)^{l}\check{f}_{j-\frac12}\right)=0,
\end{equation}
where $\check{f}_{j+\frac12}$ is a 
consistent approximation of $f(\rho)V(R)$ at interface $x_{j+1/2}$. Here, we consider again the Lax-Friedrichs numerical flux defined by
\begin{equation}\label{nlflux1}
\check{f}_{j+\frac12}=\frac12\left(f(\rho_{j+\frac12}^{h,-})V(R_{j+\frac12}^{h,-})+f(\rho_{j+\frac12}^{h,+})V(R_{j+\frac12}^{h,+})+\alpha(\rho_{j+\frac12}^{h,-}-\rho_{j+\frac12}^{h,+})  \right),
\end{equation}
with $\alpha:=\max_\rho|\partial_\rho(f(\rho)V(\rho))|$ and where $R_{j+\frac12}^{h,-}$ and $R_{j+\frac12}^{h,+}$ are the left and right approximations of $R(x,t)$ at the interface $x_{j+\frac12}$. 
Since $R$ is defined by a convolution, we naturally set
$R_{j+\frac12}^{h,-}=R_{j+\frac12}^{h,+}=R(x_{j+\frac12},t)=:R_{j+\frac12}$, so that \eqref{nlflux1} can be written as
\begin{equation}\label{nlflux2}
\check{f}_{j+\frac12}:=\check{f}(\rho_{j+\frac12}^{h,-},\rho_{j+\frac12}^{h,+},R_{j+\frac12})=
\frac12\left((f(\rho_{j+\frac12}^{h,-})+f(\rho_{j+\frac12}^{h,+}))V(R_{j+\frac12})+\alpha(\rho_{j+\frac12}^{h,-}-\rho_{j+\frac12}^{h,+})  \right).
\end{equation}
Next, we propose to approximate the integral term in (\ref{DG_nl}) using the following high-order Gauss-Legendre quadrature technique,
\begin{eqnarray}\label{nl_intR}
 \int_{I_{j}}f(\rho^h) V(R)\frac{d}{dx}v^{(l)}(\zeta_j(x))dx&=&\int_{-1}^1f(\rho^h\left(\frac{\Delta x}{2}y+x_j,t\right))V\left(R\left(\frac{\Delta x}{2}y+x_j,t\right)\right)(v^{(l)})'(y)dy \notag\\
 &=& \sum_{e=1}^{N_G}w_ef(\rho^{h} \left(\hat{x}_e,t\right))V\left(R\left(\hat{x}_e,t\right)\right)(v^{(l)})'(y_e),
\end{eqnarray}
where we have set 
$\hat{x}_e=\frac{\Delta x}{2}y_e+x_j\in I_j$, $y_e$ being the Gauss-Legendre quadrature points ensuring that the quadrature formula is exact for polynomials of order less or equal to $2N_G-1$. \\
It is important to notice that the DG formulation for the non-local conservation law \eqref{DG_nl} requires the computation of the extra integral terms $R_{j+\frac12}$ in  \eqref{nlflux2} and $R\left(\hat{x}_e,t\right)$ in \eqref{nl_intR} for each quadrature point, which increases the computational cost of the strategy. For $\rho^{h}(x,t)|_{I_j}\in P^k(I_j)$, we can compute these terms as follows for both the LWR and sedimentation non-local models considered in this paper. \\
\ \\
{\it Non-local LWR model.} For the non-local LWR model,  we impose the condition  $N\Delta x=\eta$ for some $N\in\mathbb{N}$, so that we have 
\begin{eqnarray*}
R_{j+\frac12}&=&\int_{x_{j+\frac12}}^{x_{j+\frac12}+\eta}\omega_\eta(y-{x_{j+\frac12}})\rho^h(y,t)dy=\sum_{i=1}^N\int_{I_{j+i}}
\omega_\eta(y-{x_{j+\frac12}})\rho^h_{j+i}(y,t)dy\\
&=&\sum_{i=1}^N \sum_{l=0}^{k}c^{(l)}_{j+i}(t)\int_{I_{j+i}}\omega_\eta(y-{x_{j+\frac12}})v^{(l)}(\zeta_{j+i}(y))dy\\
&=&\sum_{i=1}^N \sum_{l=0}^{k}c^{(l)}_{j+i}(t)\underbrace{\frac{\Delta x}{2}\int_{-1}^1\omega_\eta\left(\frac{\Delta x}{2}y+(i-\frac12)\Delta x\right)v^{(l)}(y)dy}_{\Gamma_{i,l}}=\sum_{i=1}^N \sum_{l=0}^{k}c^{(l)}_{j+i}\Gamma_{i,l},
\end{eqnarray*}
while for each quadrature point $\hat{x}_e$ we have
\begin{eqnarray*}
R\left(\hat{x}_e,t\right)&=&\int_{\hat{x}_e}^{\hat{x}_e+\eta}\omega_\eta(y-\hat{x}_e)\rho^h(y,t)dy= \underbrace{\int_{\hat{x}_e}^{x_{j+\frac12}}\omega_\eta(y-\hat{x}_e)\rho_{j}^h(y,t)dy}_{\Gamma_a}+\\
&&\sum_{i=1}^{N-1}\underbrace{\int_{I_{j+i}}\omega_\eta(y-\hat{x}_e)\rho_{j+i}^h(y,t)dy}_{\Gamma_b^i} +\underbrace{\int_{x_{j+N-\frac12}}^{\hat{x}_e+\eta}\omega_\eta(y-\hat{x}_e)\rho_{j+N}^h(y,t)dy}_{\Gamma_c}.
\end{eqnarray*}
The three integrals $\Gamma_a$, $\Gamma_b^i$ and $\Gamma_c$ can be computed with the same change of variable as before, namely 
\begin{eqnarray*}
\Gamma_a&=&\frac{\Delta x}{2}\int_{y_e}^1\omega_\eta\left(\frac{\Delta x}{2}(y-y_e)\right)\sum_{l=0}^kc^{(l)}_j(t)v^{(l)}(y)dy\\
&=&\sum_{l=0}^kc^{(l)}_j(t)\underbrace{\frac{\Delta x}{2}\int_{y_e}^1\omega_\eta\left(\frac{\Delta x}{2}(y-y_e)\right)v^{(l)}(y)dy}_{\Gamma_{0,l}^{(e)}}=\sum_{l=0}^kc^{(l)}_j \Gamma_{0,l}^{(e)}, \\
\Gamma_b^i&=&\frac{\Delta x}{2}\int_{-1}^1\omega_\eta\left(\frac{\Delta x}{2}(y-y_e)+i\Delta x\right)\sum_{l=0}^kc^{(l)}_{j+i}(t)v^{(l)}(y)dy\\
&=&\sum_{l=0}^kc^{(l)}_{j+i}\underbrace{\frac{\Delta x}{2}\int_{-1}^1\omega_\eta\left(\frac{\Delta x}{2}(y-y_e)+i\Delta x\right)v^{(l)}(y)dy}_{\Gamma_{i,l}^{(e)}}=\sum_{l=0}^kc^{(l)}_{j+i} \Gamma_{i,l}^{(e)}, \\
\Gamma_c&=&\frac{\Delta x}{2}\int_{-1}^{y_e}\omega_\eta\left(\frac{\Delta x}{2}(y-y_e)+N\Delta x\right)\sum_{l=0}^kc^{(l)}_{j+N}v^{(l)}(y)dy\\
&=&\sum_{l=0}^kc^{(l)}_{j+N}(t)\underbrace{\frac{\Delta x}{2}\int_{-1}^{y_e}\omega_\eta\left(\frac{\Delta x}{2}(y-y_e)+N\Delta x\right)v^{(l)}(y)dy}_{\Gamma_{N,l}^{(e)}}=\sum_{l=0}^kc^{(l)}_{j+N} \Gamma_{N,l}^{(e)}.
\end{eqnarray*}
Finally we can compute $R\left(\hat{x}_e,t\right)$ as 
$$R\left(\hat{x}_e,t\right)=\sum_{i=0}^N\sum_{l=0}^kc^{(l)}_{j+i} (t)\Gamma_{i,l}^{(e)},$$
in order to evaluate \eqref{nl_intR}. \\
\ \\
{\it Non-local sedimentation model.} 
Considering now the non-local sedimentation model, we impose $N\Delta x=2\eta$ for some $N\in\mathbb{N}$, so that we have 
$$R_{j+\frac12}= \int_{x_{j+\frac12}-2\eta}^{x_{j+\frac12}+2\eta}\omega_\eta(y-x_{j+\frac12})\rho^h(y,t)dy=\sum_{i=-N+1}^N\sum_{l=0}^kc^{(l)}_{j+i}(t)\Gamma_{i,l},$$
and for each quadrature point $\hat{x}_e$,
$$R\left(\hat{x}_e,t\right)=\int_{\hat{x}_e-2\eta}^{\hat{x}_e+2\eta}\omega_\eta(y-\hat{x}_e)\rho^h(y,t)dy=\sum_{i=-N}^N\sum_{l=0}^kc^{(l)}_{j+i}(t) \Gamma_{i,l}^{(e)},$$
with 
\begin{eqnarray*}
\Gamma_{-N,l}^{(e)}&=& \frac{\Delta x}{2}\int_{y_e}^1\omega_\eta\left(\frac{\Delta x}{2}(y-y_e)-N\Delta x\right)v^{(l)}(y)dy,\\
\Gamma_{i,l}^{(e)}&=& \frac{\Delta x}{2}\int_{-1}^1\omega_\eta\left(\frac{\Delta x}{2}(y-y_e)+i\Delta x\right)v^{(l)}(y)dy,\\
\Gamma_{N,l}^{(e)}&=& \frac{\Delta x}{2}\int_{-1}^{y_e}\omega_\eta\left(\frac{\Delta x}{2}(y-y_e)+N\Delta x\right)v^{(l)}(y)dy.
\end{eqnarray*}
\ \\
{\bf Remark.} \textcolor{black}{In order to compute integral terms in \eqref{DG_nl} as accurately as possible, the 
integrals $R_{j+\frac12}$ and $R\left(\hat{x}_e,t\right)$ above, and in particular, the coefficients $\Gamma_{i,l}$,
must be calculated exactly or using a suitable quadrature formula accurate to at least $\mathcal{O}(\Delta x ^{l+p})$ where $p$ is the degree of the convolution term $\omega_{\eta}$. 
It is important to notice that the coefficients  can be precomputed and stored in order to save computational time.}\\
\textcolor{black}{
Finally the semi-discrete scheme \eqref{DG_nl} can be discretized in time using the (TVD) third-order Runge-Kutta scheme \eqref{rk3}, under 
the CFL condition 
$$\frac{\Delta t}{\Delta x}\max_\rho|\partial_\rho (f(\rho)V(\rho))|\leq \mathrm{C}_{\mathrm{CFL}}= \frac{1}{2k+1},$$
where $k$ is the degree of the polynomial. }
%
%
\section{Construction of FV schemes for non-local conservation laws} \label{sec:fv_nl}
 
Let us now extend the FV-WENO strategy of Section \ref{sec:fvweno1} to the non-local case. 
We first integrate \eqref{nl_claw} over $I_j$ to obtain 
 \begin{equation*}\label{fv_nl_int}
\frac{d}{dt}\bar{\rho}(x_j,t)=-\frac{1}{\Delta x}\left( f(\rho(x_{j+\frac12},t))V(R(x_{j+\frac12},t))-f(\rho(x_{j-\frac12},t))V(R(x_{j-\frac12},t))\right),\quad \forall j\in\mathbb{Z},
\end{equation*}
so that the semi-discrete discretization can be written as
\begin{equation}\label{fv_semi_2}
\frac{d}{dt}\bar{\rho}_j(t)=-\frac{1}{\Delta x}\left( \check{f}_{j+\frac12}-\check{f}_{j-\frac12}\right),\qquad \forall j\in\mathbb{Z},
\end{equation}           
where the use of the Lax-Friedrichs numerical flux gives
\begin{equation*}
\check{f}_{j+\frac12}=\check{f}(\rho_{j+\frac12}^l,\rho_{j+\frac12}^r,R_{j+\frac12})=\frac12\left((f(\rho_{j+\frac12}^l)+f(\rho_{j+\frac12}^r))V(R_{j+\frac12})+\alpha(\rho_{j+\frac12}^l-\rho_{j+\frac12}^r)  \right).
\end{equation*}
Recall that $\rho_{j+\frac12}^l$ and $\rho_{j+\frac12}^r$ are the left and right WENO high-order reconstructions at point $x_{j+\frac12}$. \\
\ \\
At this stage, it is crucial to notice that in the present FV framework, the approximate solution is piecewise constant on each cell $I_j$, so that a naive evaluation of the convolution terms may lead to a loss of high-order accuracy. Let us illustrate this. Considering for instance
the non-local LWR model and using that $\rho(x,t)$ is piecewise constant on each cell naturally leads to
\begin{eqnarray*}
R_{j+\frac12}&=&\int_{x_{j+\frac12}}^{x_{j+\frac12}+\eta}\omega_\eta(y-{x_{j+\frac12}})\rho(y,t)dy=\sum_{i=1}^N\int_{I_{j+i}}
\omega_\eta(y-{x_{j+\frac12}})\rho(y,t)dy\\
&=&\Delta x \sum_{i=1}^N\bar{\rho}_{j+i}\int_{I_{j+i}}\omega_\eta(y-{x_{j+\frac12}})dy,
\end{eqnarray*}
which does not account for the high-order WENO reconstruction. In order to overcome this difficulty, 
we propose to approximate the value of $\rho(x,t)$ using quadratic polynomials in each cell.
This strategy is detailed for each model in the following two subsections.

\subsection{Computation of $R_{j+\frac12}$ for the non-local LWR model}
In order to compute the integral $$
R_{j+\frac12}=\int_{x_{j+\frac12}}^{x+\eta}\omega_\eta(y-x_{j+\frac12})\rho(y,t)dy,
$$
we propose to consider a reconstruction of $\rho(x,t)$ on $I_j$ by taking advantage of the high-order WENO reconstructions 
$\rho_{j-\frac12}^r$ and $\rho_{j+\frac12}^l$ at the boundaries of $I_j$, as well as the approximation of the cell average $\bar{\rho}_j^n$. More precisely, we propose  
to define a polynomial $P_j(x)$ of degree 2 on 
$I_j$ by 
$$
P_j(x_{j-\frac12})=\rho_{j-\frac12}^r,\quad P_j(x_{j+\frac12})=\rho_{j+\frac12}^l,\quad \frac{1}{\Delta x}\int_{I_j}P_j(x)dx=\bar{\rho}_j^n,
$$
which is very easy to handle. In particular, we have 
$$
P_j (x):=a_{j,0}+a_{j,1}\left( \frac{x-x_j}{\Delta x/2} \right)+a_{j,2}\left(3 \left( \frac{x-x_j}{\Delta x/2} \right)^2-1\right)/2,\qquad x \in I_j,
$$
with 
$$
a_{j,0}=\bar{\rho}_j^n,\qquad a_{j,1}=\frac12\left(\rho_{j+\frac12}^l-\rho_{j-\frac12}^r \right), \qquad a_{j,2}=\frac12\left(\rho_{j+\frac12}^l+\rho_{j-\frac12}^r \right)-\bar{\rho}_j^n.
$$ 
Observe that we have used the same polynomials as in the DG formulation, i.e., $P_j(x)=\sum_{l=0}^2a_{j,l}v^{(l)}(\zeta_{j}(y))$.
With this, $R_{j+\frac12}$ can be computed as 
\begin{eqnarray}\label{lwr_Rj12}
R_{j+\frac12}&=&\sum_{i=1}^N\int_{I_{j+i}}\omega_\eta(y-{x_{j+\frac12}})P_{j+i}(y)dy=\sum_{i=1}^N\int_{I_{j+i}}\omega_\eta(y-{x_{j+\frac12}})\sum_{l=0}^2a_{j+i,l}v^{(l)}(\zeta_{j+i}(y))dy \notag\\ 
&=&\sum_{i=1}^N\sum_{l=0}^2a_{j+i,l}\int_{I_{j+i}}\omega_\eta(y-{x_{j+\frac12}})v^{(l)}(\zeta_{j+i}(y))dy\notag \\ \notag
&=&\sum_{i=1}^N\sum_{l=0}^2a_{j+i,l}\underbrace{\frac{\Delta x}{2}\int_{-1}^{1}\omega_\eta\left(\frac{\Delta x}{2}y+(i-\frac12)\Delta x\right)v^{(l)}(y)dy}=\sum_{i=1}^N\sum_{l=0}^2a_{j+i,l}{\Gamma}_{i,l},
\end{eqnarray}
where the coefficients ${\Gamma}_{l,i}$ are computed exactly or using a high-order quadrature approximation. 
\subsection{Computation of $R_{j+\frac12}$ for non-local sedimentation model}
As far as the non-local sedimentation model \eqref{nl_claw}-\eqref{sednl} is concerned, we have
$$
R_{j+\frac12}=\int_{-2\eta}^{2\eta}\omega_\eta(y)\rho({x_{j+\frac12}}+y,t)dy=\int_{{x_{j+\frac12}}-2\eta}^{{x_{j+\frac12}}+2\eta}\omega_\eta(y-{x_{j+\frac12}})\rho(y,t)dy.
$$
Considering again the assumption $N\Delta x=2\eta$, we get
\begin{eqnarray*}\label{sed_Rj12}
R_{j+\frac12}&=&\sum_{i=-N+1}^N\int_{I_{j+i}}\omega_\eta(y-{x_{j+\frac12}})P_{j+i}(y)dy=\sum_{i=-N+1}^N\sum_{l=0}^2a_{j+i,l}\int_{I_{j+i}}\omega_\eta(y-{x_{j+\frac12}})v^{(l)}(\zeta_{j+i}(y))dy \notag \\
&=&\sum_{i=-N+1}^N\sum_{l=0}^2a_{j+i,l}\underbrace{\frac{\Delta x}{2}\int_{-1}^{1}\omega_\eta\left(\frac{\Delta x}{2}y+(i-\frac12)\Delta x\right)v^{(l)}(y)dy}=\sum_{i=-N+1}^N\sum_{l=0}^2a_{j+i,l}{\Gamma}_{i,l}.
\end{eqnarray*}
\\
To conclude this section, let us remark that 
the coefficients ${\Gamma}_{i,l}$ are computed for $l=0,\dots,k$ in the DG formulation, where $k$ is the degree of the polynomials in $P^k(I_j)$, while in the
FV formulation, the coefficients ${\Gamma}_{i,l}$ are computed only for $l=0,\dots,2$ due to the quadratic reconstruction of the unknown in each cell.
\textcolor{black}{
The semi-discrete scheme \eqref{fv_semi_2} is then integrated in time using the (TVD) third-order Runge-Kutta scheme \eqref{rk3}, under 
the CFL condition 
$$\frac{\Delta t}{\Delta x}\max_\rho|\partial_\rho (f(\rho)V(\rho))|\leq \mathrm{C}_{\mathrm{CFL}}<1.$$.}

%

\section{Numerical experiments} \label{sec:tests}

In this section, we propose several test cases in order to illustrate the behaviour of the RKDG and FV-WENO high-order schemes proposed in the previous Sections 
\ref{sec:dg_nl} and \ref{sec:fv_nl} for the numerical approximation of the solutions of the non-local traffic flow and sedimentation models on a bounded interval $I=[0,L]$ with boundary conditions.

We consider periodic boundary conditions for the traffic flow model, i.e. $\rho(0,t)=\rho(L,t)$ for $t\geq0$ and zero-flux boundary conditions for sedimentation model, in order to simulate a batch sedimentation process. More precisely, we assume that $\rho(x,t)=0$ for $x\leq0$ and $\rho(x,t)=1$ for $x\geq1.$ 
{
Given a uniform partition of $[0,L]$ $\{I_{j}\}_{j=1}^{M}$ with $\Delta x= L/M$, in order to compute the numerical fluxes $\check{f}_{j+1/2}$ for $j=0,\dots,M+1$, 
we define $\rho_{j}^{n}$ in the ghost cells as follow: for the traffic model, 
$$\rho_{-1}^{n}:=\rho_{M-1}^{n}, \quad \rho_{0}^{n}:=\rho_{M}^{n}, \quad \rho_{M+i}^{n}:=\rho_{i}^{n}, \qquad \text{ for  } i=1,\dots,N;$$
and for the sedimentation model 
$$\rho_{i}^{n}:=0,  \text{ for  } i=-N,\dots,0, \quad \text{ and } \quad \rho_{M+i}^{n}:=1, \qquad \text{ for  } i=1,\dots,N.$$ }

A key element of this section will be the computation of the Experimental Order of Accuracy (EOA) of the proposed strategies, which is expected to coincide with the theoretical order of accuracy given by the high-order reconstructions involved in the corresponding numerical schemes. Let us begin with a description of the practical computation of the EOA. \\ 
Regarding the RKDG schemes, if $\rho^{\Delta x}(x,T)$ and $\rho^{\Delta x/2}(x,T)\in V_h^k$  are the solutions computed with $M$ and $2M$ mesh cells respectively,  the L${}^1$-error is computed by
\begin{align*}
e(\Delta x)&=\|\rho^{\Delta x}(x,T)-\rho^{\Delta x/2}(x,T)\|_{L{}^1} \\
&=\sum_{j=1}^M \int_{I_{2j-1}} | \rho^{\Delta x}_{j}(x,T)-\rho^{\Delta x/2}_{2j-1}(x,T)|dx+
\int_{I_{2j}} | \rho^{\Delta x}_j(x,T)-\rho^{\Delta x/2}_{2j}(x,T)|dx,
\end{align*}
where the integrals are computed with a high-order quadrature formula. As far as the FV schemes are concerned, the L${}^1$-error is computed as
$$e(\Delta x)=\|\rho^{\Delta x}(x,T)-\rho^{\Delta x/2}(x,T)\|_{L{}^1}=\Delta x \sum_{j=1}^M | \rho^{\Delta x}_{j}(x_j,T)-\tilde{\rho}_{j}(x_j)|dx,$$
where $\tilde{\rho}_{j}(x)$ is a third-degree polynomial reconstruction of $\rho^{\Delta x/2}_{j}(x,T)$ at point $x_j$, i.e., 
$$\tilde{\rho}_{j}(x_{j}) =\frac{9}{16}\left( \rho^{\Delta x/2}_{2j}+\rho^{\Delta x/2}_{2j-1}\right) -\frac{1}{16}\left(\rho^{\Delta x/2}_{2j+1}+\rho^{\Delta x/2}_{2j-2} \right).
$$
In both cases, the EOA is naturally defined by 
$$\gamma(\Delta x)=\log_2\left(e(\Delta x)/e(\Delta x/2)\right).$$

In the following numerical tests, the CFL number is taken  as $\mathrm{C}_{\mathrm{CFL}}=0.2, 0.1, 0.05$ for RKDG1, RKDG2 and RKDG3  schemes respectively, and  $\mathrm{C}_{\mathrm{CFL}}=0.6$ for 
FV-WENO3 and FV-WENO5 schemes. For RKDG3 and FV-WENO5 cases $\Delta t$ is further reduced in the accuracy tests.
\subsection{Test 1a: non-local LWR model}
\begin{figure}
\centering
\begin{tabular}{cc}
(a) & (b) \\
\includegraphics[scale=0.45]{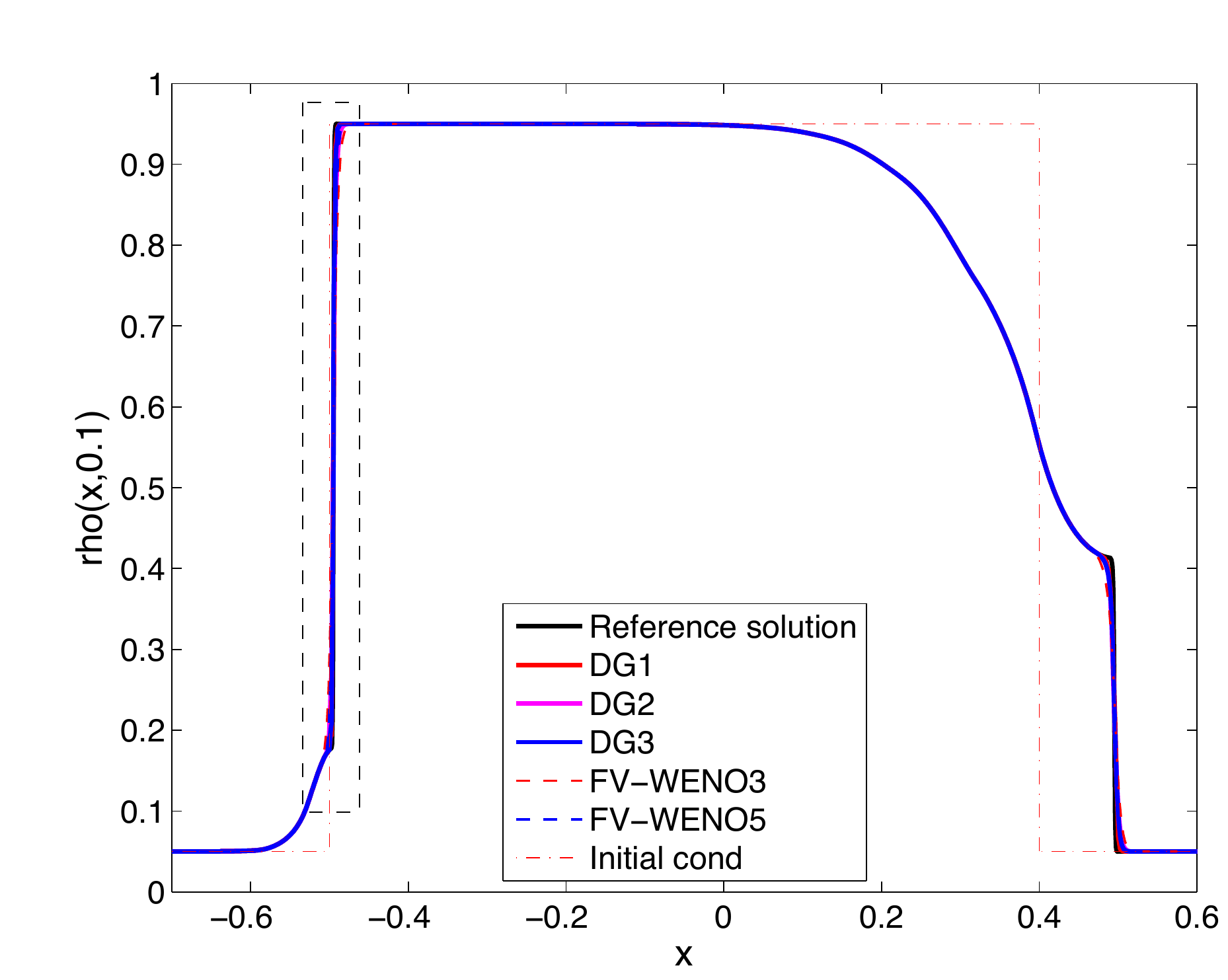}&\includegraphics[scale=0.45]{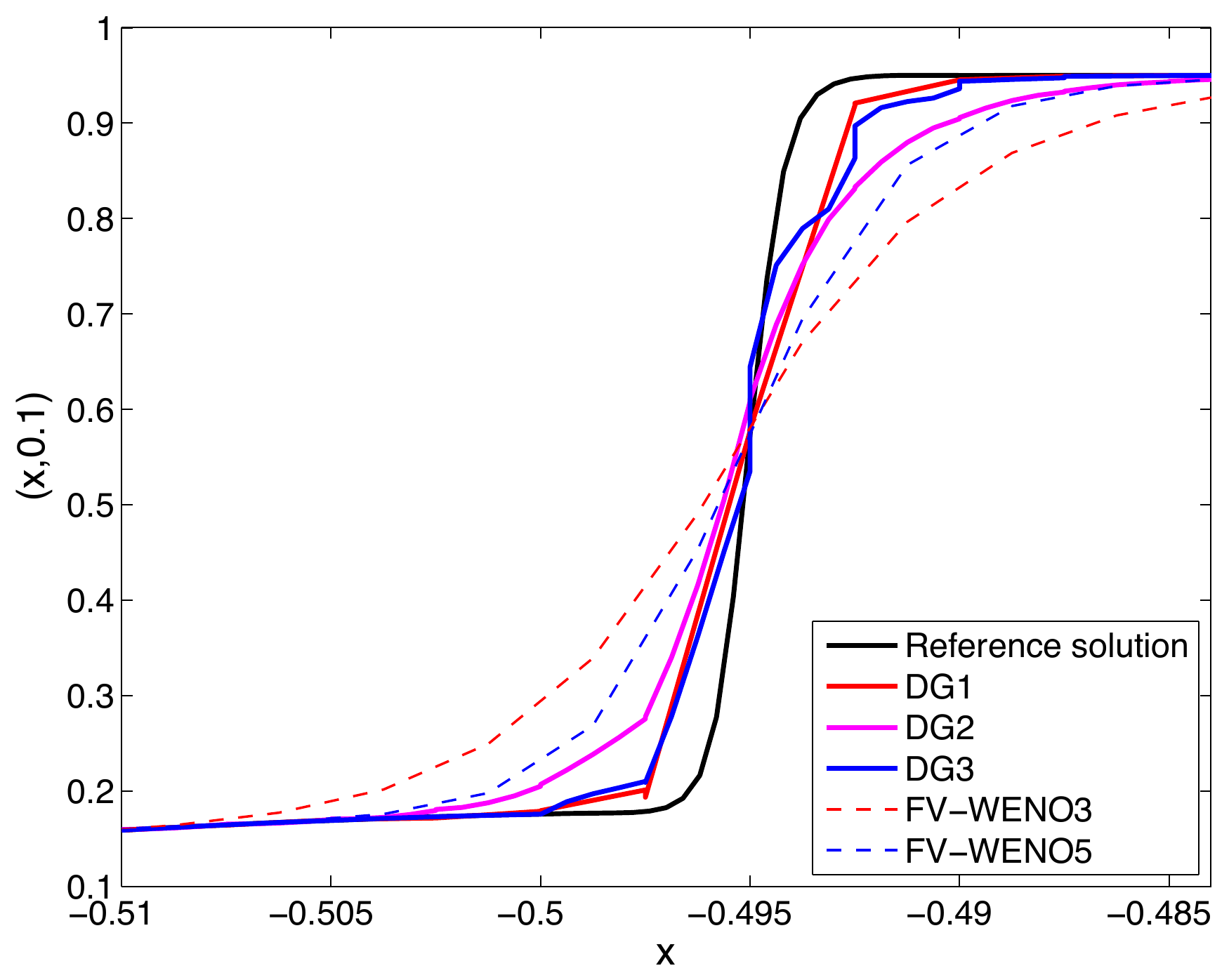}
\end{tabular}
\caption{Test 1a. (a) Solution of Riemann problem with $w_\eta(x)=3(\eta^2-x^2)/(2\eta^3)$ at $T=0.1$. We compare solutions computed with  FV-WENO and  
 RKDG schemes  using $\Delta x=1/800$. (b) Zoomed view of the left discontinuity in Fig (a). Reference solution is computed with FV-WENO5 with  $\Delta x=1/3200$. }
\label{fig_nlLWR}
\end{figure}

We consider the Riemann problem 
\begin{equation*}
\rho(x,0)= \begin{cases} 0.95 & x\in [-0.5,0.4], \\
   0.05 & \hbox{otherwise}, 
   \end{cases}
\end{equation*}
{with absorbing boundary conditions}
and compute the numerical solution of \eqref{nl_claw}-\eqref{LWRnl} at time $T=0.1$ with $\eta=0.1$ and $w_\eta(x)=3(\eta^2-x^2)/(2\eta^3)$. 
We set $\Delta x=1/800$ and compare the numerical solutions obtained with the FV-WENO3, FV-WENO5, and for  RKDG1, RKDG2  and RKDG3 we use a generalized slope limiter \eqref{slopelimiter} with $M_b=35$. 
The results  displayed in Fig. \ref{fig_nlLWR} are compared with respect to the reference solution which was obtained with FV-WENO5 scheme and $\Delta x=1/3200$. In Fig\ref{fig_nlLWR} (b) we observe that RKDG schemes are much more accurate than FV-WENO3 and even FV-WENO5.

\subsection{Test 1b: non-local LWR model}
\begin{figure}
\centering
\begin{tabular}{cc}
(a) & (b) \\
\includegraphics[scale=0.45]{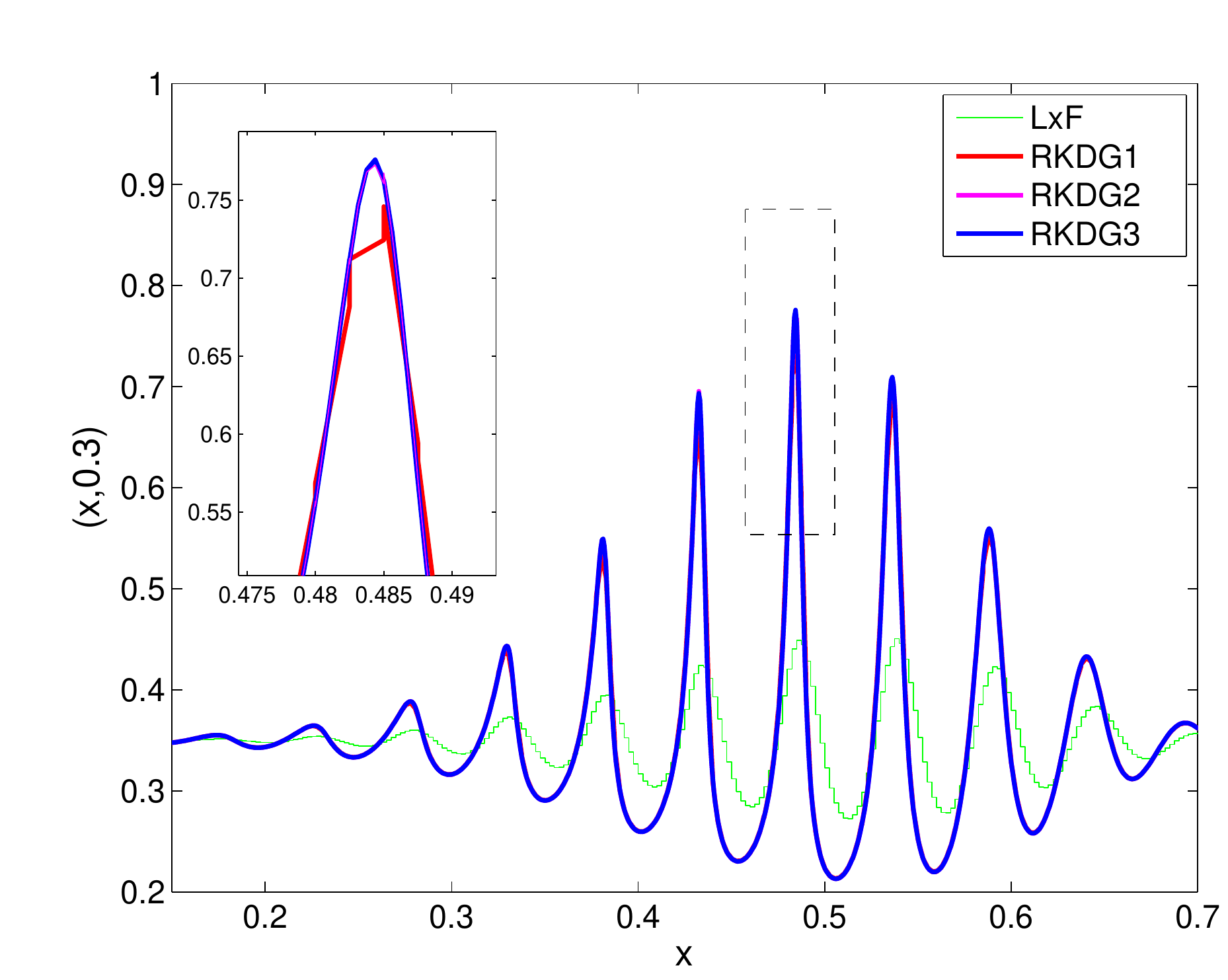}&\includegraphics[scale=0.45]{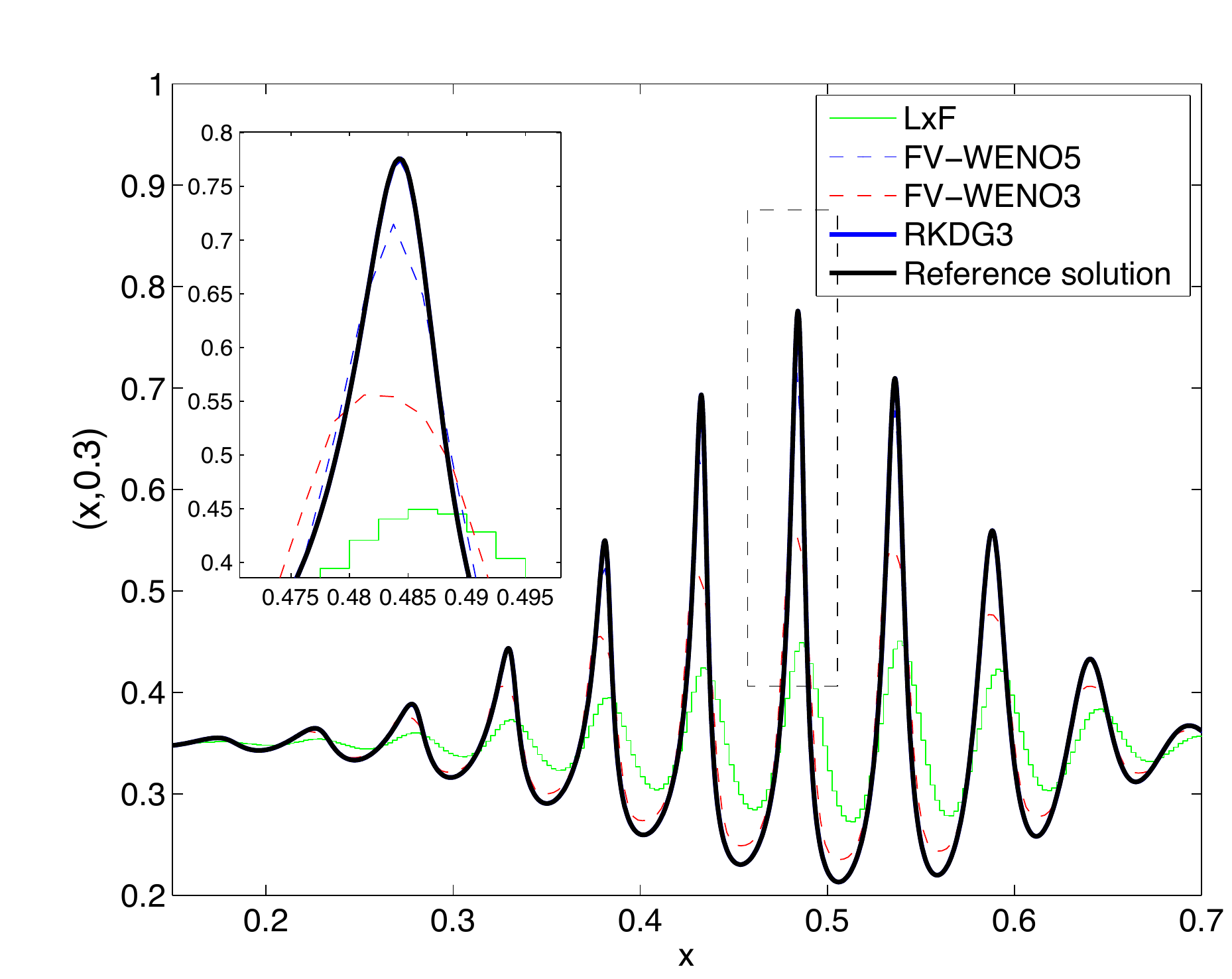}\\
\end{tabular}
\caption{ Test1b. Numerical solution non-local LWR model  with $\omega_{\eta}=2x/\eta^2$, $\eta=0.05$ and $\rho_0(x)=0.35-(x-0.5)\exp{(-2000(x-0.5)^2)}$.
We compare solutions computed with FV Lax-Friedrichs, FV-WENO  and  RKDG schemes, using $\Delta x=1/400$. Reference solution is computed with FV-WENO5 with  $\Delta x=1/3200$. }
\label{fig_nlLWR2}
\end{figure}

Now, we consider an initial condition with a small perturbation $\rho_0(x)=0.35-(x-0.5)\exp{(-2000(x-0.5)^2)}$ and an increasing kernel function 
$\omega_{\eta}=2x/\eta^2$, with $\eta=0.05$  {and periodic boundary conditions}. 
According to \cite{BG_2016, GS_2016}, the non-local LWR model is not stable with increasing kernels, in the sense that {oscillations develop in 
short time}. Fig. \ref{fig_nlLWR2} displays the numerical solution with different RKDG and FV-WENO schemes with $\Delta x=1/400$ at time $T=0.3$.
The profile provided by the first order scheme proposed in \cite{BG_2016} is also included. The reference solution is computed with a FV-WENO5 scheme with $\Delta x=1/3200$. We observe that the numerical solutions obtained with the high-order schemes provide better approximations of the oscillatory shape of the solution than the first-order scheme, and that RKDG schemes give better approximations than the FV-WENO schemes.

\subsection{Test 2: non-local sedimentation model}
\begin{figure}[t]
\begin{tabular}{cc}
\includegraphics[scale=0.45]{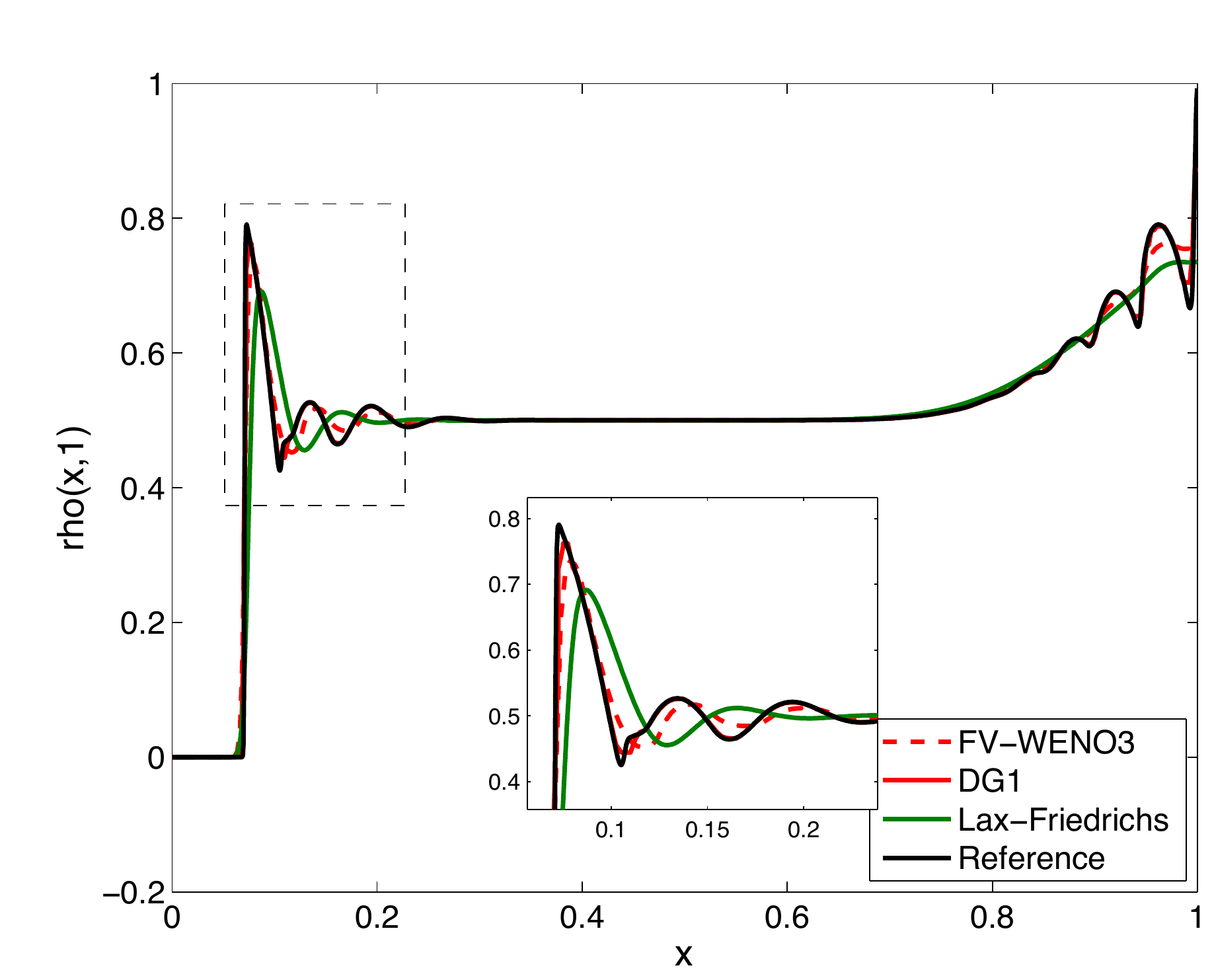}&\includegraphics[scale=0.45]{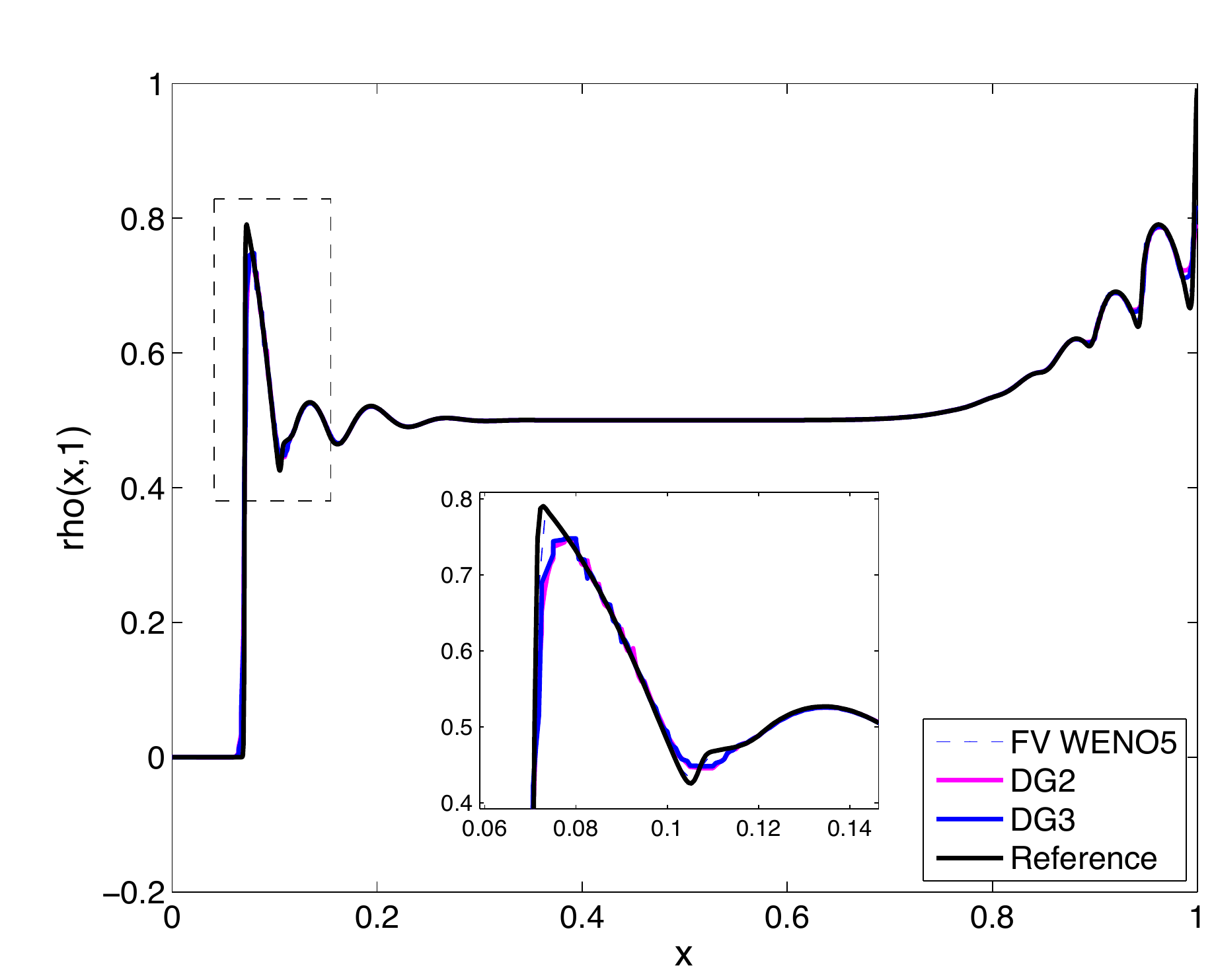} 
\end{tabular}
\caption{Test2. Riemann problem for non-local sedimentation model. Comparing numerical solution computed with FV Lax-Friedrichs, 
FV-WENO schemes and  RKDG schemes using limiter function with $M_b=50$ using $\Delta x=1/400$. Reference solution is computed with FV-WENO5 with  $\Delta x=1/3200$.}
\label{fig_nlsed}
\end{figure}

\begin{figure}
\begin{tabular}{ccc}
\fbox{\includegraphics[scale=0.3,bb=40 10 520 400]{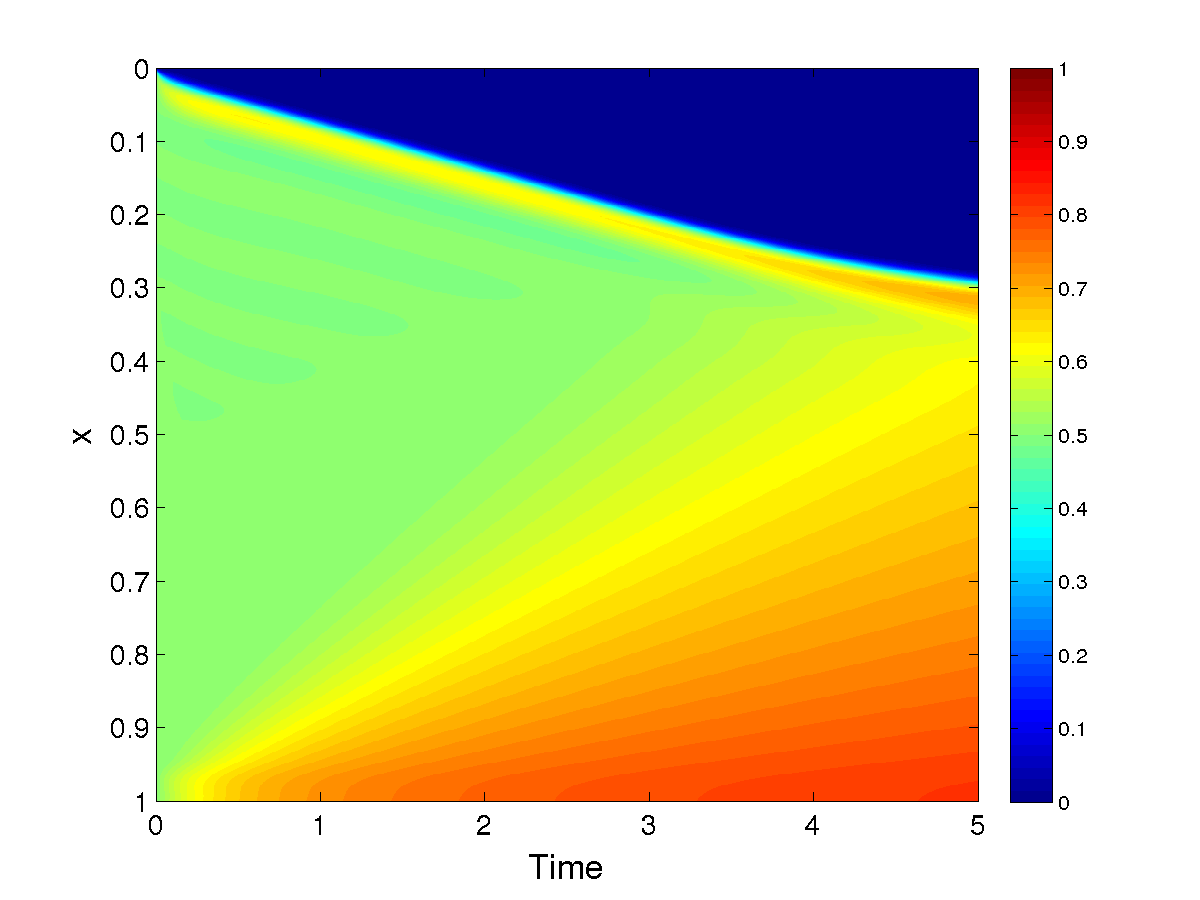}}&\fbox{\includegraphics[scale=0.3,bb=40 10 520 400]{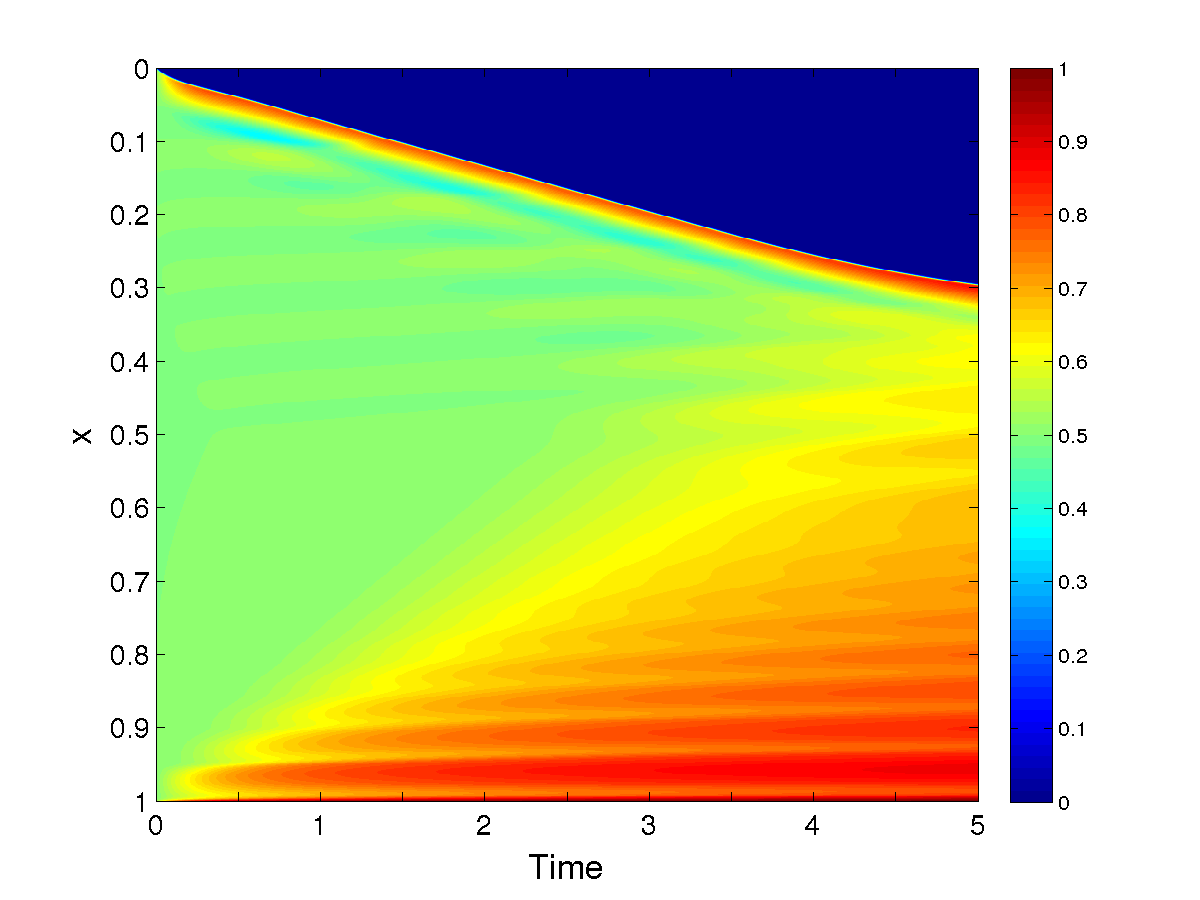}}&
\fbox{\includegraphics[scale=0.3,bb=40 10 520 400]{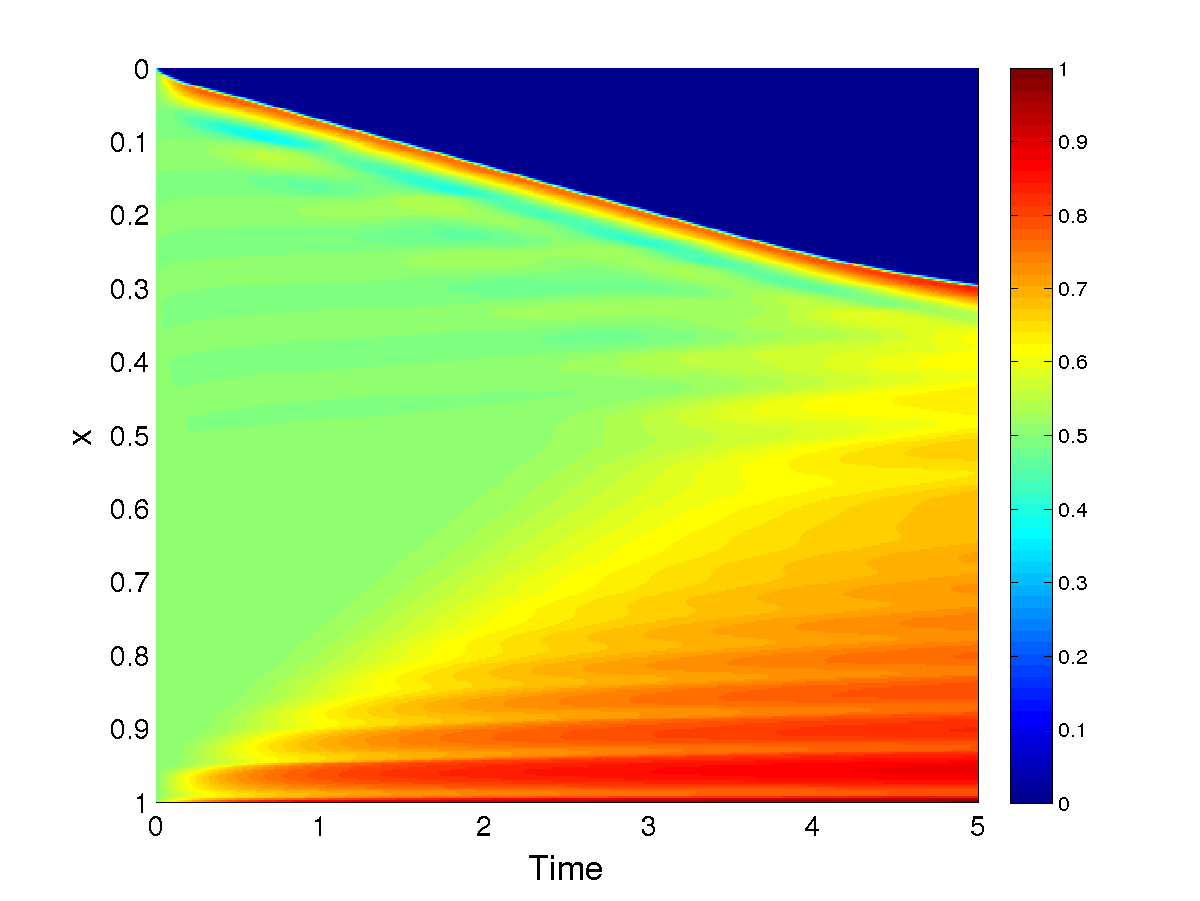}}
\end{tabular}
\caption{Test 2. Nonlocal sedimentation model. History solution for $t\in[0,5]$ computed with LxF (left) RKDG1 (center) and WENO 5 (right) with $\Delta x=1/400$.
We observe that the same kind of layering phenomenon reported in \cite{siano1979layered} seems to be observed with high order schemes RKDG1 and FV-WENO5 instead of LxF scheme.}
\label{fig_nlsedhis}
\end{figure}

We now solve \eqref{nl_claw}-\eqref{sednl}-\eqref{ker_sed} with the piece-wise constant initial datum
\begin{equation*}
\rho(x,0)= \begin{cases} 0 & x\leq 0, \\ 
0.5 & 0<x< 1, \\
   1 & x\geq1, 
   \end{cases}
\end{equation*}
{with zero-flux boundary conditions} in the interval $[0,1]$, and compute the solution at  time $T=1$, with parameters $\alpha=1$, $n=3$ and $a=0.025$. 
We set $\Delta x=1/400$ and compute the solutions with different RKDG and FV-WENO schemes, including the
first-order Lax-Friedrichs scheme used in \cite{BBKT_2011}. The results displayed in Fig. \ref{fig_nlsed} are compared to 
a reference solution computed with FV-WENO5 and  $\Delta x=1/3200$.  Compared to the reference solution, we observe that 
 RKDG1 is more accurate than FV-WENO3, and FV-WENO5 more accurate than RKDG3 and RKDG2 (Fig. \ref{fig_nlsed}).
In particular,  we observe that the numerical solutions obtained with the high-order schemes provide better approximations of the oscillatory shape of the solution than the first order scheme.
These oscillations can possibly be explained as {\it layering phenomenon} in sedimentation \cite{siano1979layered}, which denotes a traveling staircases pattern, looking 
 as several distinct bands of different concentrations. We observe in Fig. \ref{fig_nlsedhis}, where the evolution of $\rho^{h}(\cdot,t)$ is displayed for $t\in[0,5]$,
 that this layering phenomenon 
 is observed with high order schemes, in this case RKDG1 and FV-WENO5, instead of Lax-Friedrichs scheme.  
 
\subsection{Test 3: Experimental Order of Accuracy for {\it local} conservation laws} 
In this subsection we compute the EOA for local conservation laws. More precisely, we first consider the 
advection equation $\rho_t+\rho_x=0,$ and the initial datum $\rho_0(x)=0.5+0.4\sin(\pi x)$ for $x\in[-1,1]$ with periodic boundary conditions. We compare the numerical approximations 
obtained with $\Delta x=1/M$ and $M=20,40,80,160,320,640$ at time $T=1$. RKDG schemes are used without limiting function and $L^1$-errors are collected in Table \ref{EOA_advection}.
We observe that the correct EOA is obtained for the different schemes, moreover, we observe that the $L^{1}$-errors obtained with  FV-WENO3 and RKDG1 are comparable, and also the 
$L^{1}$-errors  for  FV-WENO5 and RKDG3.
We also consider the nonlinear local LWR model   $\rho_t+(\rho(1-\rho))_x=0,$ with initial datum $\rho_0(x)=0.5+0.4\sin(\pi x)$ for $x\in[-1,1]$ and periodic boundary conditions. A shock wave appears at time $T=0.3$. We compare the numerical approximations obtained with $\Delta x=1/M$ and $M=20,40,80,160,320,640$ at  time $T=0.15$. RKDG schemes are used without limiters and $L^1$-errors are given in Table \ref{EOA_LWR}. We again  observe that the correct EOA is obtained for the different schemes, and the $L^{1}$-errors obtained with  FV-WENO3 and RKDG1 are comparable, and also the 
$L^{1}$-errors  for  FV-WENO5 and RKDG3.

\begin{table}[h!]
  \centering
  \begin{tabular}{ |c|c|c|c|c|c|c||c|c|c|c|}\hline
     &  \multicolumn{2}{|c|}{$k=1$}  &  \multicolumn{2}{|c|}{$k=2$}  &  \multicolumn{2}{|c||}{$k=3$} &  \multicolumn{2}{|c|}{FV-WENO$3$}&  \multicolumn{2}{|c|}{FV-WENO$5$}   \\ \hline
    $\tiny{1/\Delta x}$ & $L^1$-err  &\tiny{$\gamma(\Delta x)$} & $L^1$-err  &\tiny{$\gamma(\Delta x)$}   & $L^1$-err  &\tiny{$\gamma(\Delta x)$} & $L^1$-err  &\tiny{$\gamma(\Delta x)$} 
    & $L^1$-err  &\tiny{$\gamma(\Delta x)$} \\ \hline\hline
    $20$   &1.80E-03   &--       &9.70E-05  &  --      & 1.00E-06 & --     &3.72E-02&--&4.63E-04&-- \\ 
    $40$   &4.30E-04   &{2.06}&1.28E-05   & 2.92  & 6.04E-08&  4.05& 1.17E-02&1.66 & 1.41E-05  &  5.03  \\ 
    $80$   &1.06E-04   & 2.02 & 1.63E-06  & 2.96 &  3.92E-09& 3.94 & 2.85E-03& 2.04&  4.33E-07&  5.02 \\ 
    $160$&2.64E-05    & 2.00 & 2.06E-07  & 2.98 &  2.58E-10&  3.92& 4.24E-04& 2.74& 1.35E-08 &  5.00 \\
    $320$ &6.61E-06  & 2.00 &2.59E-08    & 2.99&1.51E-11&3.32&3.06E-05  & 3.79& 1.15E-11& 5.20\\ \hline
  \end{tabular}
  \caption{Test 3: Advection equation.  Experimental order of accuracy for  FV-WENO schemes and RKDG schemes without Limiter 
  at $T=1$  with initial condition $\rho_0(x)=0.5+0.4\sin(\pi x)$ on $x\in[-1,1]$.   }
  \label{EOA_advection}
\end{table} 

\begin{table}[h]
  \centering
  \begin{tabular}{ |c|c|c|c|c|c|c||c|c|c|c|}\hline
     &  \multicolumn{2}{|c|}{$k=1$}  &  \multicolumn{2}{|c|}{$k=2$}  &  \multicolumn{2}{|c||}{$k=3$} &  \multicolumn{2}{|c|}{WENO$3$}&  \multicolumn{2}{|c|}{WENO$5$}   \\ \hline
    $\tiny{1/\Delta x}$ & $L^1$-err  &\tiny{$\gamma(\Delta x)$} & $L^1$-err  &\tiny{$\gamma(\Delta x)$}   & $L^1$-err  &\tiny{$\gamma(\Delta x)$} & $L^1$-err  &\tiny{$\gamma(\Delta x)$} 
    & $L^1$-err  &\tiny{$\gamma(\Delta x)$} \\ \hline\hline
    $20$   &2.03E-03   & -- & 1.88E-04 &-- &5.03E-06  &--&6.92E-03&--&2.80E-04&-- \\ 
    $40$    & 4.98E-04 & 2.02 &  3.06E-05& 2.62 &  2.81E-07 &  4.16& 2.37E-03&1.77&1.41E-05 &  4.31       \\ 
    $80$     & 1.23E-04 & 2.01  &  4.71E-06& 2.69 & 1.72E-08  & 4.03 & 5.95E-04&1.99 &6.14E-07&4.52   \\ 
     $160$  &3.08E-05  &2.00   & 7.11E-07 &   2.72 &  1.05E-09& 4.02& 8.20E-05& 2.86&  2.62E-08& 4.55 \\ 
     $320$  &7.71E-06 &2.00  &1.06E-07&2.74&6.62E-11&3.99 & 5.47E-06&3.90 & 8.61E-10&  4.92  \\\hline
  \end{tabular}
  \caption{Test 3: Classical LWR equation. Experimental order of accuracy  for  FV-WENO schemes and RKDG schemes without Limiter 
  at $T=0.15$  with initial condition $\rho_0(x)=0.5+0.4\sin(\pi x)$ on $x\in[-1,1]$.  }
  \label{EOA_LWR}
\end{table}


\subsection{Test 4: Experimental Order of Accuracy for the {\it non-local} problems}
We now consider the non-local LWR and sedimentation models. Considering the non-local LWR model, we compute the solution of \eqref{nl_claw}-\eqref{LWRnl} with initial data $\rho_0(x)=0.5+0.4\sin(\pi x)$ for $x\in[-1,1]$, { periodic boundary conditions}, with $\eta=0.1$ and $\Delta x=1/M$ and $M=80,160,320,640,1280$ at final time $T=0.15$.  The results are given for different kernel functions $w_\eta(x)$ in Table \ref{EOA_nl_LWR}. For RKDG schemes, we obtain the correct EOA. For the FV-WENO schemes, the EOA is also correct thanks to the in-cell quadratic reconstructions used to 
compute the non-local terms. \\
For the non-local sedimentation model, we compute the solution of \eqref{nl_claw}-\eqref{sednl}-\eqref{ker_sed} with initial data $\rho_0(x)=0.8\sin(\pi x)^{10}$ for $x\in[0,1]$ with $\eta=0.025$, $\alpha=1$, $V(\rho)=(1-\rho)^3$ and $\Delta x=1/M$ with $M=100,200,400,800,1600,3200$ at final time $T=0.04$. The results are given in Table \ref{EOA_nl_LWR}.

 \begin{table}
  \centering
  \begin{tabular}{ |c||c|c|c|c|c|c|c| }\hline
    & &  \multicolumn{2}{|c|}{$w_\eta(x)=1/\eta$}  &  \multicolumn{2}{|c|}{$w_\eta(x)=2(\eta-x)/\eta^2$}  &  \multicolumn{2}{|c|}{$w_\eta(x)=3(\eta^2-x^2)/(2\eta^3)$}  \\ \hline
    &$1/\Delta x$ & $L^1-$error  &$\gamma(\Delta x)$ & $L^1-$error  &$\gamma(\Delta x)$  &  $L^1-$error  &$\gamma(\Delta x)$ \\ \hline\hline
     \multirow{5}{*}{\rotatebox[origin=c]{90}{ FVWENO3}}&$40$  & 1.49E-03   & --   & 1.64 -03   & --    & 1.78E-03     & -- \\ 
    &$80$  & 3.27E-04   &   2.18  &3.39E-04    & 2.27     &  3.53E-04     &   2.33   \\ 
    &$160$ & 4.64E-05   & 2.81   & 4.63E-05   & 2.87     &4.69E-05       &      2.91 \\ 
    &$320$ &  3.37E-06    &  3.78   &  3.44E-06   &  3.74     &    3.57E-06     &    3.71 \\
     &$640$ & 2.29E-07    & 3.87   & 2.41E-07   & 3.83     &  2.50E-07      & 3.83     \\ \hline
     
    \multirow{5}{*}{\rotatebox[origin=c]{90}{FVWENO5}}&$40$  &7.54E-06  & --  &  8.65E-06 & -- & 9.69E-06  & \\ 
   & $80$  & 3.36E-07     &4.48    &3.75E-07  & 4.52 & 4.17E-07    & 4.53\\ 
   & $160$ &1.21E-08     &4.79   &1.44E-08   &4.70  & 1.64E-08   & 4.66\\ 
   & $320$ &2.80E-10     & 5.43  & 4.52E-10 & 4.99 &   5.85E-10 & 4.81 \\
    & $640$ &6.71E-12    &  5.38   &1.85E-11  & 4.61 &  2.70E-11   & 4.43  \\ \hline\hline

      \multirow{5}{*}{\rotatebox[origin=c]{90}{RKDG $k=1$}}&$40$  & 4.72E-04  &  --& 4.75E-04   & --& 4.76E-04   & --\\ 
   & $80$  &1.17E-04      &   2.00    & 1.18E-04  & 1.99  &  1.19E-04 &1.99  \\ 
   & $160$ & 2.94E-05    &   2.00   & 2.97E-05  & 1.99  & 2.99E-05   &1.99  \\  
   & $320$ & 7.36E-06    &   1.99  & 7.45E-06  & 1.99 &   7.49E-06  & 1.99 \\ 
    & $640$ & 1.84E-06     &  1.99    & 1.86E-06  & 1.99 & 1.87E-06    &  1.99  \\ \hline  
 
        \multirow{5}{*}{\rotatebox[origin=c]{90}{RKDG $k=2$}}&$40$  & 1.81E-05    &  --& 1.21E-05  & --& 1.15E-05    & --\\ 
   & $80$  &  2.63E-06     & 2.78 & 1.77E-06  & 2.77 &  1.65E-06 & 2.80\\ 
   & $160$ &  3.70E-07  & 2.82  &  2.54E-07 & 2.79  & 2.40E-07  & 2.78 \\  
   & $320$ &  3.96E-08    & 2.89  & 3.58E-08   & 2.82  & 3.43E-08   & 2.80\\ 
    & $640$ & 5.20E-09  & 2.93  & 4.03E-09 & 2.86 & 4.81E-09  &2.83   \\ \hline       
     
         \multirow{5}{*}{\rotatebox[origin=c]{90}{RKDG $k=3$}}&$40$  &1.93E-07   &  --&1.86E-07   & --& 1.68E-07    & --\\ 
   & $80$  &1.21E-08     & 3.99 &1.20E-08   &3.95  &1.05E-08   & 4.00\\ 
   & $160$ &7.55E-10   & 4.00 &7.65E-10   &3.97   &6.65E-10   & 3.98\\  
   & $320$ &4.87E-11   & 3.95 & 4.97E-11  & 3.94 & 4.35E-11  &3.93 \\ 
    & $640$ & 3.41E-12 &  3.83 &3.51E-12  & 3.82  & 3.15E-12  & 3.78  \\ \hline      
    
  \end{tabular}
  \caption{Test 4: Non-local LWR model. Experimental order of accuracy. Initial condition $\rho_0(x)=0.5+0.4\sin(\pi x)$, $\eta=0.1$, numerical solution at $T=0.15$ for
  FV-WENO and RKDG schemes without generalized slope limiter.}
  \label{EOA_nl_LWR}
\end{table}

\begin{table}[h]
  \centering
  \begin{tabular}{ |c|c|c|c|c|c|c||c|c|c|c|}\hline
     &  \multicolumn{2}{|c|}{$k=1$}  &  \multicolumn{2}{|c|}{$k=2$}  &  \multicolumn{2}{|c||}{$k=3$} &  \multicolumn{2}{|c|}{WENO$3$}&  \multicolumn{2}{|c|}{WENO$5$}   \\ \hline
    $\tiny{1/\Delta x}$ & $L^1$-err  &\tiny{$\gamma(\Delta x)$} & $L^1$-err  &\tiny{$\gamma(\Delta x)$}   & $L^1$-err  &\tiny{$\gamma(\Delta x)$} & $L^1$-err  &\tiny{$\gamma(\Delta x)$} 
    & $L^1$-err  &\tiny{$\gamma(\Delta x)$} \\ \hline\hline
    $100$   & 7.98E-05 & --     & 2.90E-06 & --&4.25E-07   &--  & 2.72E-04 & -- & 2.65E-06 &  --  \\   
    $200$   &1.98E-05  &2.0    & 4.64E-07 & 2.6  & 3.95E-08  &3.5  &6.02E-05 &2.1  & 1.01E-07 &   4.7     \\ 
    $400$   & 4.95E-06&  2.0   & 7.75E-08 & 2.6  & 3.37E-09  &3.6  & 7.87E-06 &2.9  & 3.22E-09 &   4.9 \\ 
     $800$  & 1.23E-06 &  2.0  & 1.23E-08 &  2.6 &  2.47E-10  &3.8  &5.50E-07 &3.8  &1.11E-10  & 4.8 \\ 
     $1600$& 3.09E-07 &  2.0  & 2.02E-09 & 2.6  &  1.22E-11 & 4.2 & 4.64E-08 &3.7   &4.02E-12  & 4.7 \\\hline
  \end{tabular}
  \caption{Test 4: Non-local sedimentation problem. Initial condition $\rho_0(x)=0.8\sin(\pi x)^{10}$. 
  Experimental order of accuracy at $T=0.04$ with $f(\rho)=\rho(1-\rho)$, $V(\rho)=(1-\rho)^3$, $\omega_\eta(x):=\eta^{-1}K(\eta^{-1}x)$ with $\eta=0.025$.}
  \label{EOA_nl_sed}
\end{table}

\section*{Conclusion}
In this paper we developed high order numerical approximations of the solutions of non-local conservation laws in one space dimension motivated by application to traffic flow and sedimentation models.
We propose to design Discontinuous Galerkin (DG) schemes which can be applied in a natural way and Finite Volume WENO (FV- WENO) schemes where we have used quadratic polynomial reconstruction in each cell to evaluate the convolution term in order to obtain the high-order accuracy. The numerical solutions obtained with high-order schemes provide better approximations of the oscillatory shape of the solutions 
than the first order schemes. We also remark that  DG schemes are more accurate but expensive, while FV-WENO are less accurate but also less expensive. The present work thus establishes the preliminary basis for a deeper study and optimum programming of the methods, so that efficiently plots could then be considered. 
\section*{Acknowledgements}
This research was conducted while
LMV was visiting Inria Sophia Antipolis M\'editerran\'ee and the Laboratoire de Math\'ematiques de Versailles, 
with the support of  Fondecyt-Chile project  11140708 and  ``Poste Rouge'' 2016 of  CNRS -France. 

The authors thank R\'egis Duvigneau for useful discussions and advices.

{ \small
 \bibliography{CGV}

\def\ocirc#1{\ifmmode\setbox0=\hbox{$#1$}\dimen0=\ht0 \advance\dimen0
  by1pt\rlap{\hbox to\wd0{\hss\raise\dimen0
  \hbox{\hskip.2em$\scriptscriptstyle\circ$}\hss}}#1\else {\accent"17 #1}\fi}
\begin{thebibliography}{10}

\bibitem{abramowitz1966handbook}
M.~Abramowitz, I.~A. Stegun, et~al.
\newblock Handbook of mathematical functions.
\newblock {\em Applied mathematics series}, 55:62, 1966.

\bibitem{ACG2015}
A.~Aggarwal, R.~M. Colombo, and P.~Goatin.
\newblock Nonlocal systems of conservation laws in several space dimensions.
\newblock {\em SIAM J. Numer. Anal.}, 53(2):963--983, 2015.

\bibitem{AHP2016}
D.~Amadori, S.-Y. Ha, and J.~Park.
\newblock On the global well-posedness of \{BV\} weak solutions to the
  kuramoto-sakaguchi equation.
\newblock {\em Journal of Differential Equations}, 262(2):978 -- 1022, 2017.

\bibitem{AmadoriShen2012}
D.~Amadori and W.~Shen.
\newblock An integro-differential conservation law arising in a model of
  granular flow.
\newblock {\em J. Hyperbolic Differ. Equ.}, 9(1):105--131, 2012.

\bibitem{AmbrosioGangbo2008}
L.~Ambrosio and W.~Gangbo.
\newblock Hamiltonian {ODE}s in the {W}asserstein space of probability
  measures.
\newblock {\em Comm. Pure Appl. Math.}, 61(1):18--53, 2008.

\bibitem{Amorim2012}
P.~Amorim.
\newblock On a nonlocal hyperbolic conservation law arising from a gradient
  constraint problem.
\newblock {\em Bull. Braz. Math. Soc. (N.S.)}, 43(4):599--614, 2012.

\bibitem{Amorim2015}
P.~Amorim, R.~Colombo, and A.~Teixeira.
\newblock On the numerical integration of scalar nonlocal conservation laws.
\newblock {\em ESAIM M2AN}, 49(1):19--37, 2015.

\bibitem{BBKT_2011}
F.~Betancourt, R.~B{\"u}rger, K.~H. Karlsen, and E.~M. Tory.
\newblock On nonlocal conservation laws modelling sedimentation.
\newblock {\em Nonlinearity}, 24(3):855--885, 2011.

\bibitem{BG_2016}
S.~Blandin and P.~Goatin.
\newblock Well-posedness of a conservation law with non-local flux arising in
  traffic flow modeling.
\newblock {\em Numer. Math.}, 132(2):217--241, 2016.

\bibitem{Carrillo2016}
J.~A. Carrillo, S.~Martin, and M.-T. Wolfram.
\newblock An improved version of the {H}ughes model for pedestrian flow.
\newblock {\em Math. Models Methods Appl. Sci.}, 26(4):671--697, 2016.

\bibitem{CS_2001}
B.~Cockburn and C.-W. Shu.
\newblock Runge-{K}utta discontinuous {G}alerkin methods for
  convection-dominated problems.
\newblock {\em J. Sci. Comput.}, 16(3):173--261, 2001.

\bibitem{ColomboGaravelloMercier2012}
R.~M. Colombo, M.~Garavello, and M.~L\'ecureux-Mercier.
\newblock A class of nonlocal models for pedestrian traffic.
\newblock {\em Mathematical Models and Methods in Applied Sciences},
  22(04):1150023, 2012.

\bibitem{ColomboHertyMercier2011}
R.~M. Colombo, M.~Herty, and M.~Mercier.
\newblock Control of the continuity equation with a non local flow.
\newblock {\em ESAIM Control Optim. Calc. Var.}, 17(2):353--379, 2011.

\bibitem{ColomboMercier2012}
R.~M. Colombo and M.~L{\'e}cureux-Mercier.
\newblock Nonlocal crowd dynamics models for several populations.
\newblock {\em Acta Math. Sci. Ser. B Engl. Ed.}, 32(1):177--196, 2012.

\bibitem{ColomboMarcellini2015}
R.~M. Colombo and F.~Marcellini.
\newblock Nonlocal systems of balance laws in several space dimensions with
  applications to laser technology.
\newblock {\em J. Differential Equations}, 259(11):6749--6773, 2015.

\bibitem{CrippaMercier2012}
G.~Crippa and M.~L\'ecureux-Mercier.
\newblock Existence and uniqueness of measure solutions for a system of
  continuity equations with non-local flow.
\newblock {\em Nonlinear Differential Equations and Applications NoDEA}, pages
  1--15, 2012.

\bibitem{GS_2016}
P.~Goatin and S.~Scialanga.
\newblock Well-posedness and finite volume approximations of the {LWR} traffic
  flow model with non-local velocity.
\newblock {\em Netw. Heterog. Media}, 11(1):107--121, 2016.

\bibitem{Gottlich2014}
S.~G{\"o}ttlich, S.~Hoher, P.~Schindler, V.~Schleper, and A.~Verl.
\newblock Modeling, simulation and validation of material flow on conveyor
  belts.
\newblock {\em Applied Mathematical Modelling}, 38(13):3295 -- 3313, 2014.

\bibitem{GottlichSchindler2015}
S.~G{\"o}ttlich and P.~Schindler.
\newblock Discontinuous {G}alerkin {M}ethod for {M}aterial {F}low {P}roblems.
\newblock {\em Math. Probl. Eng.}, pages Art. ID 341893, 15, 2015.

\bibitem{gottlieb2009high}
S.~Gottlieb, D.~I. Ketcheson, and C.-W. Shu.
\newblock High order strong stability preserving time discretizations.
\newblock {\em Journal of Scientific Computing}, 38(3):251--289, 2009.

\bibitem{Keimer2015}
M.~Gugat, A.~Keimer, G.~Leugering, and Z.~Wang.
\newblock Analysis of a system of nonlocal conservation laws for
  multi-commodity flow on networks.
\newblock {\em Netw. Heterog. Media}, 10(4):749--785, 2015.

\bibitem{KurganovPolizzi2009}
A.~Kurganov and A.~Polizzi.
\newblock Non-oscillatory central schemes for a traffic flow model with
  {A}rrehenius look-ahead dynamics.
\newblock {\em Netw. Heterog. Media}, 4(3):431--451, 2009.

\bibitem{Perthame_book2007}
B.~Perthame.
\newblock {\em Transport equations in biology}.
\newblock Frontiers in Mathematics. Birkh\"auser Verlag, Basel, 2007.

\bibitem{PiccoliRossi2013}
B.~Piccoli and F.~Rossi.
\newblock Transport equation with nonlocal velocity in {W}asserstein spaces:
  convergence of numerical schemes.
\newblock {\em Acta Appl. Math.}, 124:73--105, 2013.

\bibitem{PiccoliTosin2011}
B.~Piccoli and A.~Tosin.
\newblock Time-evolving measures and macroscopic modeling of pedestrian flow.
\newblock {\em Arch. Ration. Mech. Anal.}, 199(3):707--738, 2011.

\bibitem{QS_2005}
J.~Qiu and C.-W. Shu.
\newblock Runge-{K}utta discontinuous {G}alerkin method using {WENO} limiters.
\newblock {\em SIAM J. Sci. Comput.}, 26(3):907--929, 2005.

\bibitem{Shu_1998}
C.-W. Shu.
\newblock Essentially non-oscillatory and weighted essentially non-oscillatory
  schemes for hyperbolic conservation laws.
\newblock In {\em Advanced numerical approximation of nonlinear hyperbolic
  equations ({C}etraro, 1997)}, volume 1697 of {\em Lecture Notes in Math.},
  pages 325--432. Springer, Berlin, 1998.

\bibitem{shu1988efficient}
C.-W. Shu and S.~Osher.
\newblock Efficient implementation of essentially non-oscillatory
  shock-capturing schemes.
\newblock {\em Journal of Computational Physics}, 77(2):439--471, 1988.

\bibitem{siano1979layered}
D.~B. Siano.
\newblock Layered sedimentation in suspensions of monodisperse spherical
  colloidal particles.
\newblock {\em Journal of Colloid and Interface Science}, 68(1):111--127, 1979.

\end{thebibliography}
\bibliographystyle{abbrv}
}

\end{document}